\documentclass[10pt]{article} % figure generation depends on this font size
\usepackage{geometry}
\usepackage{graphicx} % Required for inserting images
\usepackage{amsmath, amssymb, amsthm, amsfonts}
\usepackage{float}
\usepackage{xcolor}
\usepackage{todonotes}
\usepackage{subcaption}
\usepackage{hyperref}
\usepackage{url}
\usepackage{mwe}

\theoremstyle{remark}
\newtheorem{remark}{Remark}

\title{Searching for $J$-holomorphic curves via machine: first steps}
\author{James Rowan, Yuan Yao}
\date{\today}

\begin{document}

\maketitle

\begin{abstract}
We assemble numerical algorithms to search for $J$-holomorphic curves in symplectic manifolds. 
Each algorithm employs several different numerical techniques, each technique addressing a different aspect of the geometric problem. 
We separately consider both classical Fourier expansion and deep neural networks in our algorithms and compare their performance. 
Our algorithms take as input a smooth curve in a given homology class and search for a $J$-holomorphic curve in the same homology class.
We first verify we can produce explicitly known holomorphic curves in complex manifolds, for example the Weierstrass $\wp$ function on the torus and curves in $S^2\times S^2$ with the standard complex structure. Then we search for $J$-holomorphic curves in $S^2\times S^2$ with non-integrable almost complex structures: essentially we start with a known holomorphic curve in an integrable almost complex structure $J_0$, deform $J_0$ to a nearby nonintegrable almost complex structure $J_\epsilon$, and use our methods to find the nearby $J_\epsilon$-holomorphic curve. 
\end{abstract}

\section{Introduction}
A symplectic manifold is a manifold $M$ with a closed nondegenerate 2-form $\omega$. They originally arise as the phase space in classical mechanics \cite{arnold}, for example in the study of celestial mechanics. They also appear naturally in high energy physics \cite{mirror_symmetry} through the lens of mirror symmetry in string theory.\footnote{More than half a century of development means the references in symplectic geometry are too numerous to list. Where possible, we shall refer the reader to canonical textbooks in the subject and the references they contain.}

There was a huge boom in the development and understanding of symplectic manifolds in the 1980s after Gromov's introduction of $J$- holomorphic curves \cite{Gromov}. 
A tame almost complex structure on a symplectic manifold $(M,\omega)$ is an endomorphism $J$ of $TM$ satisfying
\begin{enumerate}
\item $J^2=-1$
\item $\omega(v,Jv)>0$ for all $v\in TM$.
\end{enumerate}
An almost complex structure $J$ is called integrable if it induces a complex structure on the symplectic manifold, otherwise it is called nonintegrable.
A symplectic manifold has an infinite dimensional contractible space of tame almost complex structures. By contrast, the space of integrable complex structures is finite dimensional. In practice this means compared to complex structures, we have a lot of freedom to choose an almost complex structure to help answer a specific question in symplectic geometry, and having a contractible space means whichever specific almost complex structure we choose will result in identical (enumerative, homological, categorical) invariants that we extract using this almost complex structure. 
We refer the reader to the canonical textbook \cite{intro_sympl} for the theory of symplectic manifolds.

Studying the space of $J$-holomorphic curves is the method by which we extract information about the symplectic manifold using an almost complex structure $J$. A $J$-holomorphic curve is a map $u:(\Sigma,j)\rightarrow (M,\omega, J)$, here $\Sigma$ is a Riemann surface\footnote{For us we will always take the Riemann surface to be closed. For the purpose of this paper when discussing $J$-holomorphic curves we do not quotient by the automorphisms of the domain. See Remark \ref{rmk:holmaps} for a detailed discussion.} with complex structure $j$, and $u$ satisfies the Cauchy-Riemann equation 
\[
J\circ du =du \circ j.
\]

We refer the reader to the canonical textbook \cite{McDuff_Salamon} and the references therein for a systematic treatise on $J$-holomorphic curves. Even when the almost complex structure is not integrable, $J$-holomorphic curves share many common properties with holomorphic curves in complex geometry: they satisfy analytic continuation, have finitely many singularities, and in nice cases come in smooth families. What distinguishes $J$-holomorphic curves from their algebraic counterparts is the flexibility of $J$. Many problems in symplectic geometry require us to choose a non-integrable $J$ and understand $J$-holomorphic curves for that particular $J$. Indeed a central paradigm in symplectic geometry in the last 40 years is given a problem (for example finding symplectic embeddings or understanding dynamics of vector fields on symplectic manifolds), first establish how to relate the problem to finding $J$-holomorphic curves for a particular kind of $J$, then work extremely hard to establish that the desired $J$-holomorphic curves exist. See Chapter 9 of \cite{McDuff_Salamon} for a nice exposition.

Finding $J$-holomorphic curves in a symplectic manifold for a general $J$ is a notoriously hard problem and is in general completely intractable. In the case where $J$ is integrable one can sometimes write down the curves by hand, however as soon as one leaves the land of integrable $J$ there is no general method by which we can write down a solution to the nonlinear Cauchy-Riemann equation by hand. 
Symplectic geometers have developed indirect methods to construct $J$-holomorphic curves, for example first finding some known holomorphic curves in algebraic settings and showing they persist as we deform the complex structure to be nonintegrable. Even then, these methods are highly implicit, producing statements along the lines of ``there exists a $J$-holomorphic curve in a given homology class'' and giving little concrete information about what the curve looks like.

There has been many instances of using numerical methods to understand geometric problems, see for example the survey on numerical methods in minimal surfaces \cite{min_surface_survey} and the references therein.
The recent introduction of physics-informed neural networks~\cite{raissi2019PINNs} has also led a surge of interest in studying geometric PDEs numerically; see for example~\cite{Pinn_ricci, 2026usersguidepinnsgeometric, 2026minimalsurfacesknotsneural, corts2026machinelearningapproachnirenberg, fang2021,douglas22a, MR4608987,MR4761208,Gerdes_2023}.

In this work we investigate using numerical methods to search for $J$-holomorphic curves in symplectic manifolds. We consider both classical numerical techniques and neural network techniques.

\subsection{Summary of results}

We focus on two classes of problems. First, as a proof of concept, we work in a symplectic manifold with an integrable $J$ in which the solutions to the $J$-holomorphic curve equation are explicitly known. We start with an initial smooth guess curve in the same homology class which is very far from being holomorphic, then run our numerical algorithms to upgrade the smooth curve to the known holomorphic curve. 

The first specific example in this class we consider is degree 2 holomorphic curves $T^2\rightarrow S^2$ that satisfy 4 marked point constraints in Section \ref{sec:wp}. Here the complex structure on $T^2$ is fixed, and we do not quotient by the automorphisms of the domain, see Remark \ref{rmk:holmaps}. With these specifications, there are two distinct holomorphic curves $\wp^{(1)}$ and $\wp^{(2)}$ that satisfy the marked point constraints, both coming from domain reparametrization/Mobius transformations of the Weierstrass $\wp$ function. We start with 16 different random smooth functions $T^2\rightarrow S^2$ and in Sections~\ref{sss:paradigm1} and~\ref{sss:paradigm2} apply our algorithms to search for the desired holomorphic curves. Out of these initial smooth functions, approximately 1/3 converge to $\wp^{(1)}$, 1/3 converge to $\wp^{(2)}$ and the remaining 1/3 do not converge. For a nice visualization of this process, see Figures~\ref{fig:p2_wp1} and~\ref{fig:p2_wp2}. This result is compelling, because a question of great interest to symplectic geometers is not just establishing existence of $J$-holomorphic curves, but also enumerating $J$-holomorphic curves (such enumerations are then assembled into symplectic invariants, like Gromov-Witten invariants or Floer homology). This result suggests our methods have the potential to answer enumerative questions in symplectic geometry.

The second specific example is the symplectic manifold $S^2\times S^2$ with the standard product complex structure. We look for holomorphic maps $S^2\rightarrow S^2\times S^2$ of bidegree $(5,3)$ that have 18 marked points. After choosing the marked points correctly, the only holomorphic curve that satisfies the marked point constraints is $z\rightarrow (z^5,z^3)$. In Section \ref{sec:5-3} we use our algorithms to recover this map starting from a smooth map of the same bidegree that does not satisfy the marked point constraints. A good visualization is in Figure \ref{fig:5-3}.

\begin{figure}[H]
    \centering
    \includegraphics[width=0.8\linewidth]{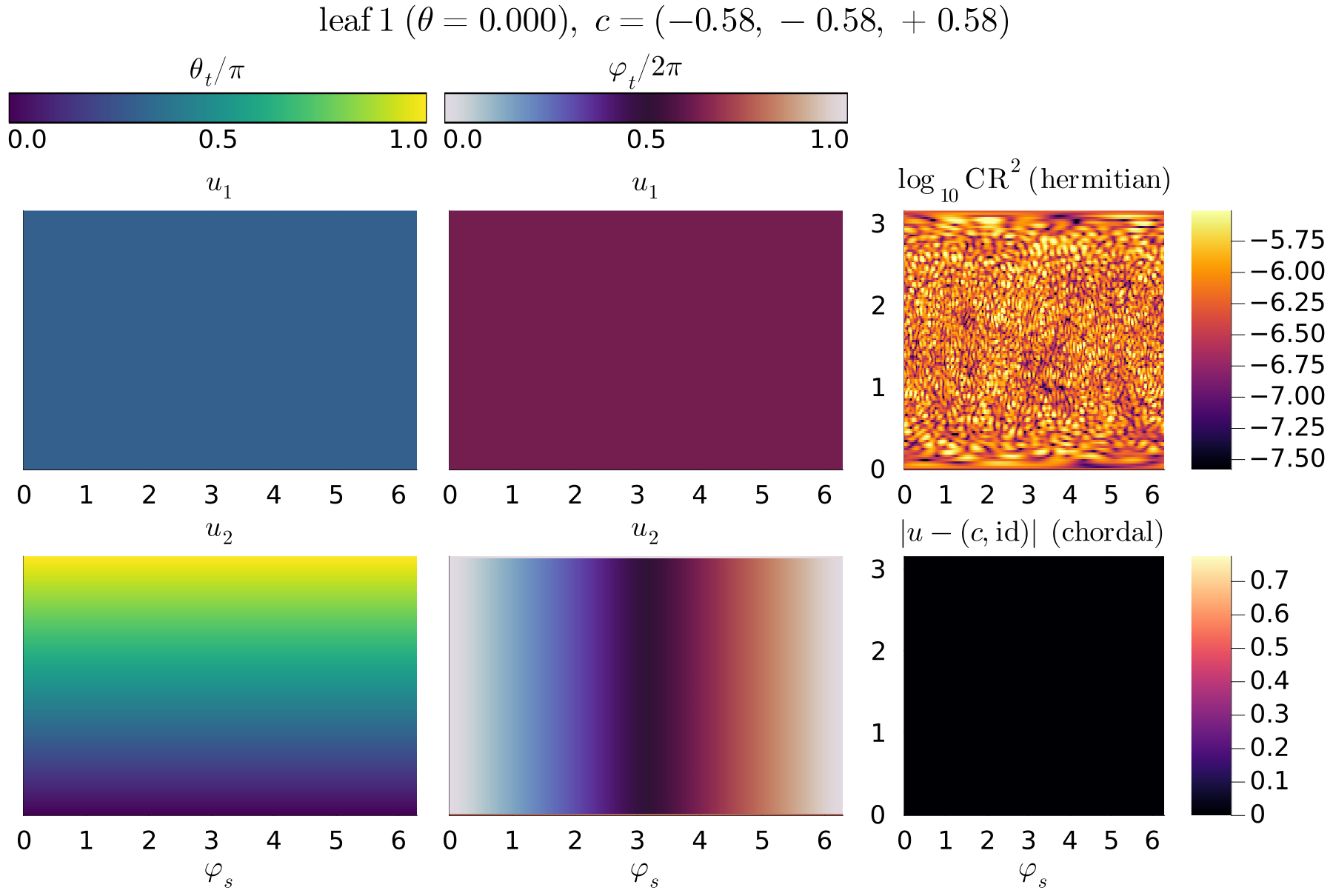}
    \caption{The initial leaf of the foliation with the degree $0$ factor's marked point equal to $c_0$; this corresponds to the holomorphic curve $u^0(z)=(c_0,z)$}
    \label{fig:01foliation_1}
\end{figure}

\begin{figure}[H]
    \centering
    \includegraphics[width=0.8\linewidth]{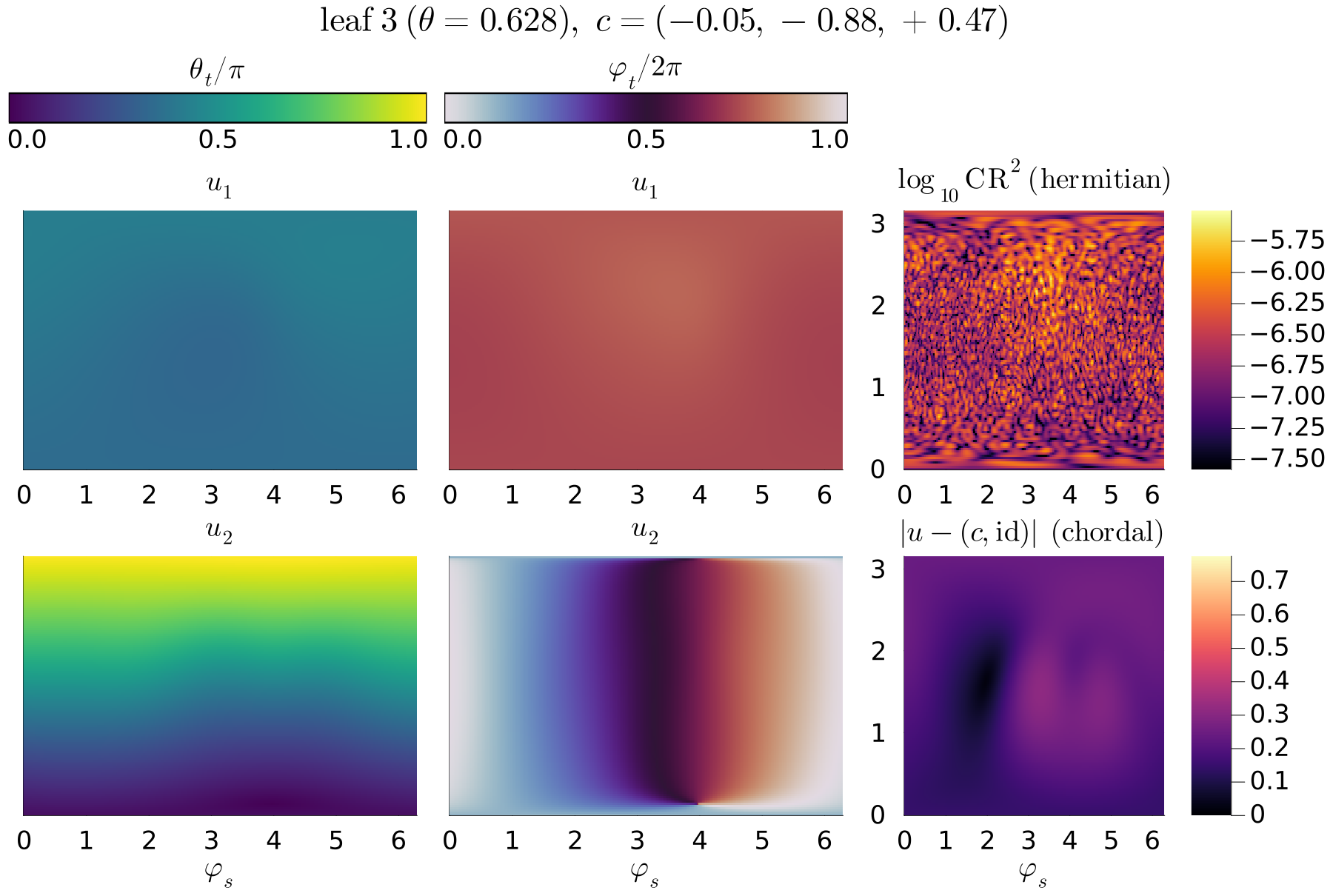}
    \caption{The leaf of the foliation with the degree $0$ factor's marked point equal to $c_2$, which is at an angle of $\frac{\pi}{5}$ along a great circle starting at $c_0$.
    Because $J$ is non-integrable, distortion away from $(c_2,z)$ is observed.}
    \label{fig:01foliation_3}
\end{figure}

\begin{figure}[H]
    \centering
    \includegraphics[width=0.8\linewidth]{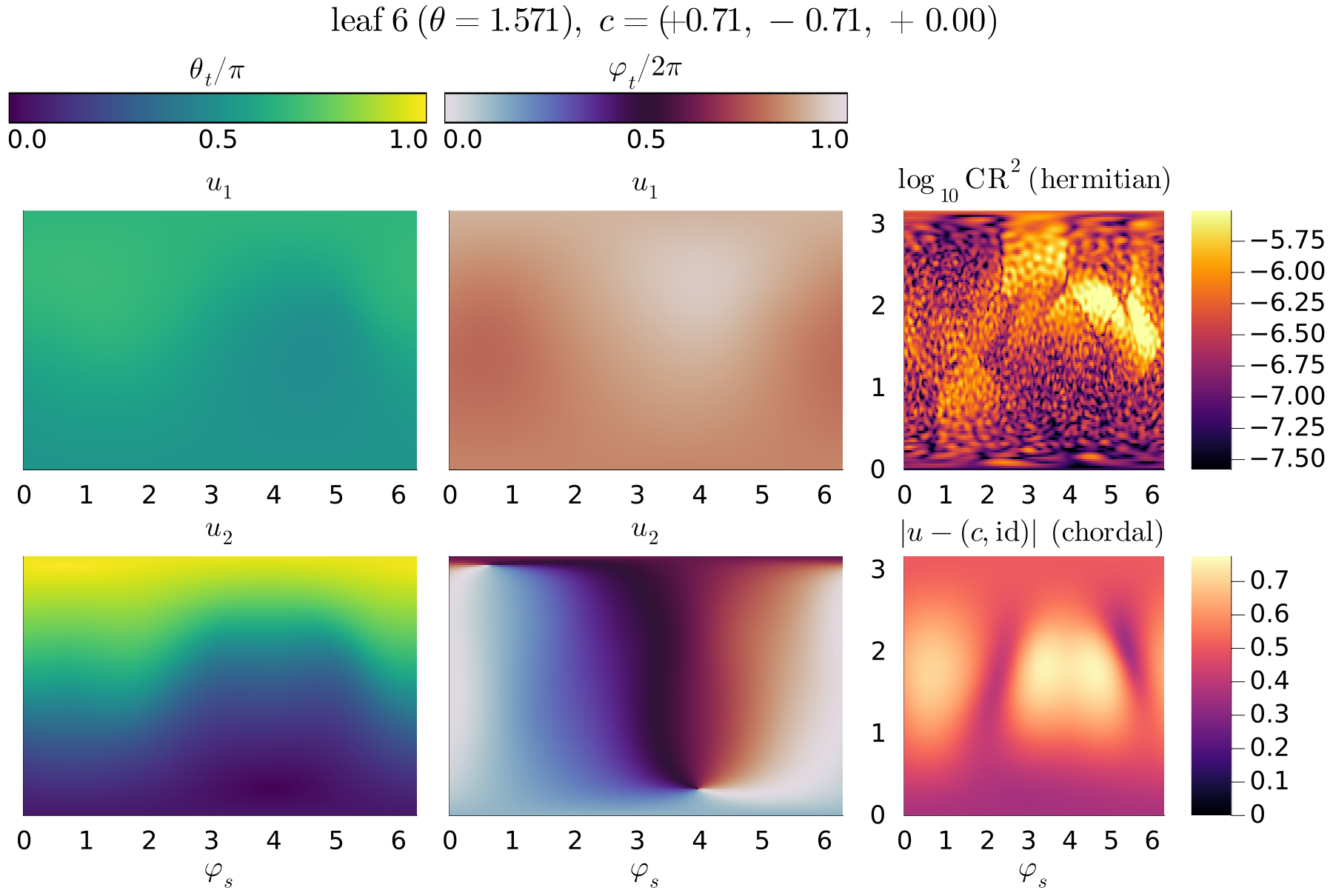}
    \caption{The leaf of the foliation with the degree $0$ factor's marked point equal to $c_6$.
    $c_6$ is at an angle of $\frac{\pi}{2}$ along a great circle starting at $c_0$, and represents the maximum amount of distortion away from the split integrable complex structure $J_0$ attainable with our $J$.}
    \label{fig:01foliation_6}
\end{figure}

\begin{figure}[H]
    \centering
    \includegraphics[width=0.8\textwidth]{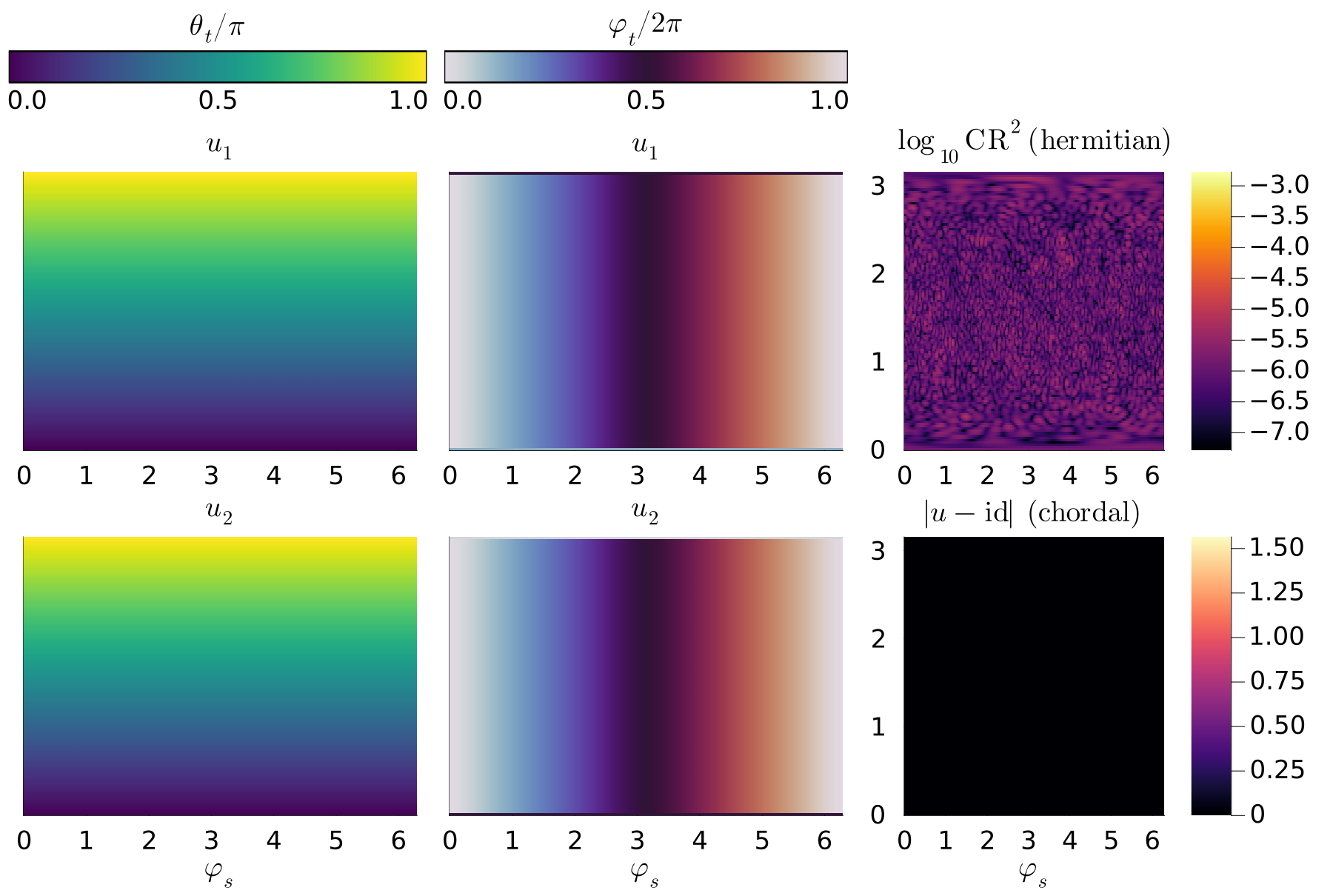}
    \caption{The learned identity map $C_0$}
    \label{fig:continuation_identity}
\end{figure}

\begin{figure}[H]
    \centering
    \includegraphics[width=0.8\textwidth]{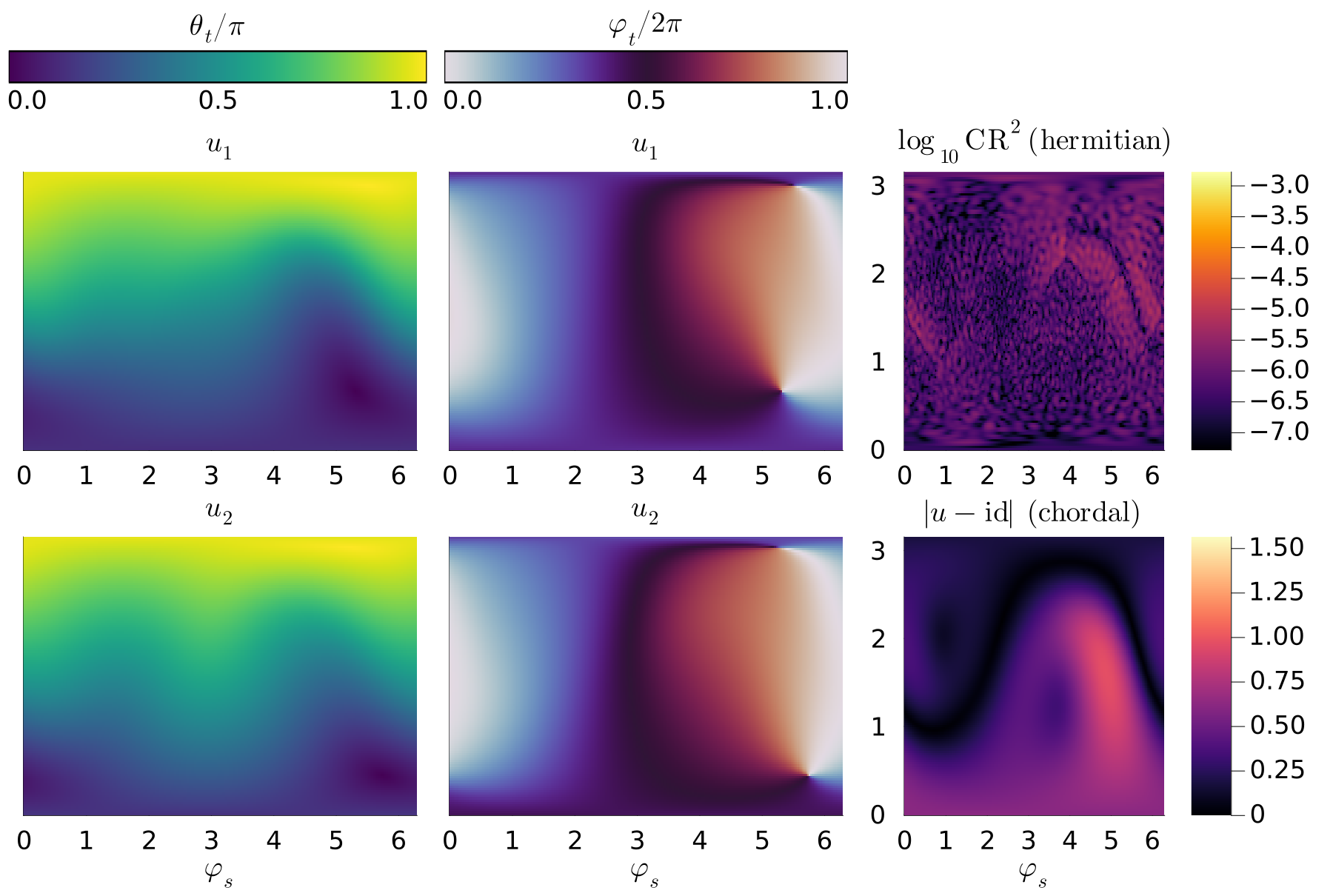}
    \caption{The curve $C_2$ with $\epsilon_2=0.6$}
    \label{fig:continuation_0.6}
\end{figure}

\begin{figure}[H]
    \centering
    \includegraphics[width=0.8\textwidth]{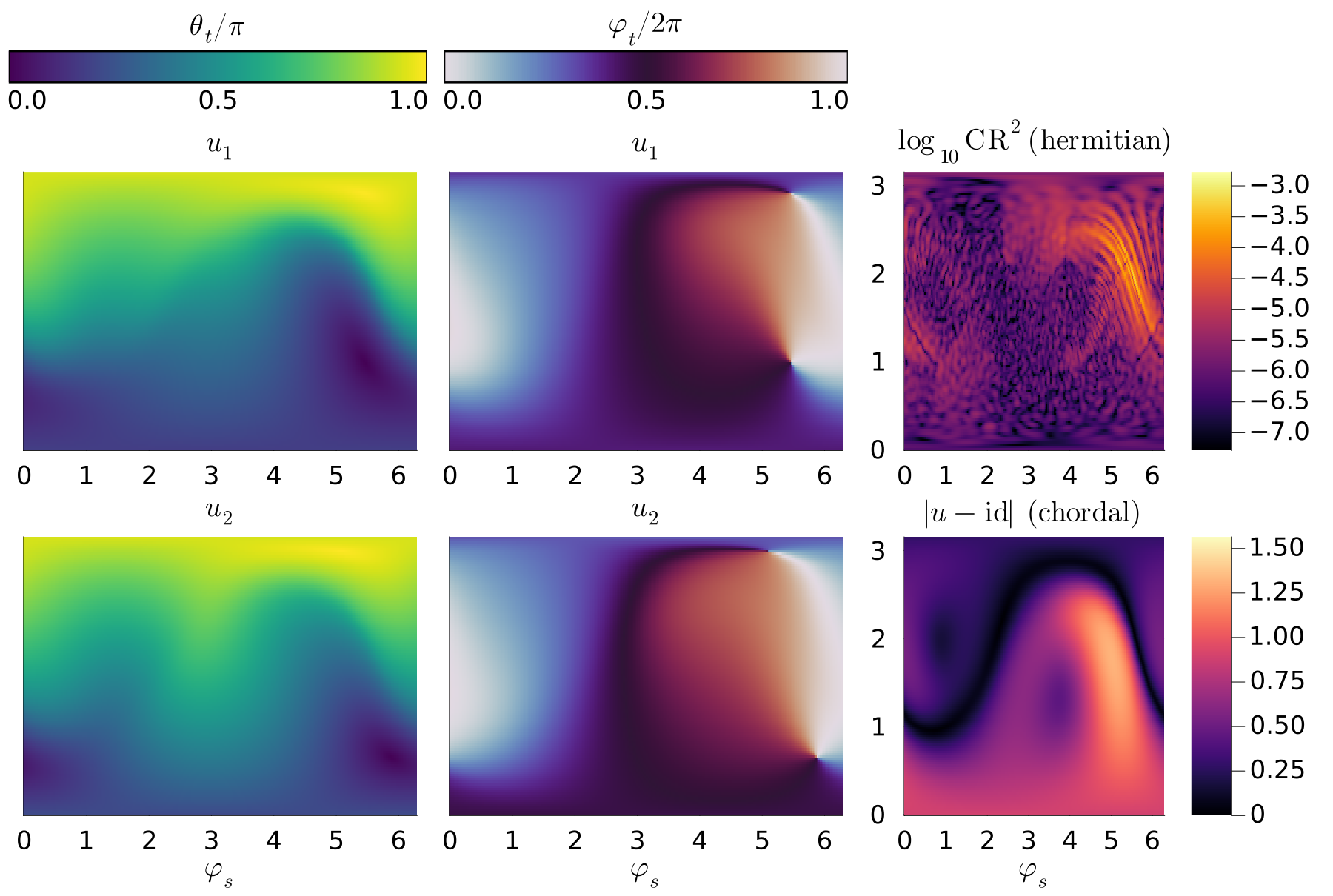}
    \caption{The curve $C_3$ with $\epsilon_3=0.9$}
    \label{fig:continuation_0.9}
\end{figure}

Once the viability of our methods has been established in these test cases, we use them to search for $J$-holomorphic curves in non-integrable settings. We study $J$-holomorphic curves in $S^2\times S^2$ of bidegree $(0,1)$ in Section \ref{ss:gromovFoliation} and bidegree $(1,1)$ in Section \ref{ss:bifurcation}. While the implementation details differ, the basic premise of both of these examples is the same: start with a $J_0$-holomorphic curve with an integrable $J_0$, deform $J_0$ to be nonintegrable $J_\epsilon$ and search for the corresponding $J_\epsilon$-holomorphic curve. For curves of degree $(1,1)$ this is a precise description of what we do in Section \ref{ss:bifurcation}, and we produce visualizations of this process in Figures~\ref{fig:continuation_identity},~\ref{fig:continuation_0.6}, and~\ref{fig:continuation_0.9}. For the $(0,1)$ curves the precise implementation is slightly different and described in detail in Section \ref{ss:gromovFoliation}, because we want to take advantage of the fact that $(0,1)$ curves form a foliation of $S^2\times S^2$. See Figures~\ref{fig:01foliation_1},~\ref{fig:01foliation_3}, and~\ref{fig:01foliation_6} for a  visualization of the results.

\subsection{Outline of methodology}

There are several distinct algorithms we use to achieve the above. Each algorithm consists of two important components.
The first component is the choice of representation for smooth functions from the domain to the target, using either Fourier modes or a deep neural network.
The specifics of these representations are described in Section \ref{sec:param}.
The next component is the method of how to go from smooth functions to $J$-holomorphic curves.
This part specifies which kind of smooth functions we use as input and the optimizer used to find the $J$-holomorphic curve, drawing on various numerical tools such as Newton iteration and stochastic gradient descent.
Our choices of how to assemble these tools is described in Section \ref{sec:training}; we will refer to this overall strategy as the paradigm. Each paradigm is more than just an arbitrary combination of methods; different paradigms represent different geometric processes by which we use to search for $J$-holomorphic curves. Schematic diagrams of these paradigms we use are given in Figures~\ref{f:paradigm1} and~\ref{f:paradigm2}. 

A question of great current discussion within the PDEs community is when and how neural networks outperform classical numerical methods.
To our knowledge, there were simply no previous builds of classical numerical methods to study $J$-holomorphic curves in the literature.
So we built both classical methods and neural network methods in this paper. The most careful comparison of pros and cons of different algorithms and paradigms is performed for the case of the Weierstrass function $\wp$ in Section \ref{sec:wp}; we give a brief description of some relevant considerations below. For the other $J$-holomorphic curves we consider, as this is a proof of concept paper, we illustrate how one particular choice of algorithm produces desired results.

\begin{remark}
There are several competing interests to balance in selecting a numerical approach.
Quasi-Newton methods can converge rapidly to high precision with a good initial guess, but their running time scales poorly in the number of model parameters.
Furthermore, adaptive regularization terms (such as the term in~\eqref{e:gradregularization} penalizing relative spikes in the gradient) interact poorly with quasi-Newton methods.
As a result, we experimented with both Physics-informed neural networks (PINNs) and traditional numerical methods.
The added representational capacity from using deep neural networks instead of a simple sum of Fourier modes has to be weighed against the slower training process and need for optimization algorithms that scale better with parameter count.
Some preliminary analysis is shown in Table~\ref{tab:hyperparametersweep} in Section~\ref{sss:hyperparametersweep}; holding optimizer and training schedule constant, a deep neural network outperformed a shallow neural network or purely linear model.
However, on the multi-initialization runs of Sections~\ref{sss:paradigm1} and~\ref{sss:paradigm2}, we had success with both PINNs and classical methods, albeit they are run under different paradigms with the particular choice of paradigm based on the hyperparameter sweep in Section ~\ref{sss:hyperparametersweep}. 
Neither clearly outperformed the other at converging to a solution, though we note that classical methods were about an order of magnitude faster and had lower final error when they did converge.
Further systematic analysis and hyperparameter tuning of both methods is warranted, in particular to better understand the difficulties in the loss landscape posed by the combination of       analytic, geometric, and topological constraints.
\end{remark}

We fully acknowledge our tool can only provide numerical solutions to the Cauchy-Riemann equation. We envision our tool to be useful to symplectic geometers in the following way: while the analytic properties of the solutions of the Cauchy-Riemann equation are in many cases well understood, the key difficulty that remains in symplectic geometry is often in constructing and explicitly enumerating the solutions of the Cauchy-Riemann equation. We envision our tool can help symplectic geometers decide with a reasonable degree of certainty whether there is a $J$-holomorphic curve in a certain symplectic manifold or not (and in nice cases give an answer to how many), and if so leverage the more traditional tools in symplectic geometry to rigorously construct these curves.

Our code is written in the Julia language and is available on Github\footnote{\url{https://github.com/jrowan/NumericalJHolomorphicCurves}}.

\subsection{List of conventions}
We clarify the conventions we used in our code. We identify $T^2=\mathbb{R}^2/\mathbb{Z}^2$ with the almost complex structure $j= \begin{pmatrix}
0 & -1\\
1& 0
\end{pmatrix}.
$
We view $S^2\subset \mathbb{R}^3$ as the unit sphere. We shall refer to the point $(0,0,1)$ as the north pole. We view the standard complex plane $\mathbb{C}\subset S^2$ with stereographic projection as $F:S^2\setminus (0,0,1)\rightarrow \mathbb{C}$ given by 
\[
F(X,Y,Z) = \left (\frac{X}{1-Z},\frac{Y}{1-Z} \right)
\]
and the complex structure on $S^2$ comes from pulling back the standard almost complex structure on $\mathbb{C}$. As such the complex structure on $S^2$ can be written using the cross product. If $v\in T_xS^2$, then $Jv = v\times x$. This is tamed by the following symplectic form: for $v_1,v_2\in T_xS^2$, $\omega(v_1,v_2) = \langle x, -v_1\times v_2\rangle$ (note this is the opposite convention as in example 4.1.5 in \cite{intro_sympl}).

When we speak of a holomorphic map $\mathbb{C}\rightarrow \mathbb{C}$, such as $z\rightarrow z^3$, we can extend both the domain and target to view it as a map $S^2\rightarrow S^2$. Similarly a meromorphic function on $T^2$ can be viewed as a holomorphic function $T^2\rightarrow S^2$.

\subsection{Concurrent Work}

We understand potentially many people are currently experimenting with machine learning and numerical methods in symplectic geometry. If you also have a project studying aspects of $J$-holomorphic curves from a numerical point of view, please feel free to send us an email and we would be happy to acknowledge you and discuss.

When our work neared the final stages of its preparation, we were delighted to discover Elliot Kienzle has independently and simultaneously built a PINN based solver for the Floer equation in Hamiltonian Floer homology. Their results are made publicly available on their website \cite{elliot}. We refer the reader to the textbook \cite{audin_damian} and the references therein for the foundations of Hamiltonian Floer homology. Viewed in the correct way, a solution to the Floer equation is a $J$-holomorphic curve whose domain is the (non-compact) infinite cylinder $\mathbb{R}\times S^1$ (we list this and $J$-holomorphic curves with more complicated domains as directions for future work in Section \ref{sec:future}). It has been long recognized that the Hamiltonian Floer chain complex, which is built out of counts of solutions to the Floer equation, contains a lot more interesting information than the Floer homology itself (the homology is already known to be isomorphic to singular homology). However, to this day we only have limited understanding of the chain complex itself because of the difficulty of directly finding solutions to Floer's equations. As such, we expect E. Kienzle's work to be of great interest to the symplectic geometry community as well.

\section{Parameterizing a smooth function} \label{sec:param}

In this section we will describe how to parametrize a smooth function $$u:X\rightarrow Y $$ with given domain and target. 
In this paper, the domain $X$ will always be $S^2$ or $T^2$. 
We will think of $T^2=[0,1]^2$ with the boundary identified. 
We will think of the domain $S^2$ as the unit sphere in $\mathbb{R}^3$. Our target $Y$ will always be $S^2$ or $S^2\times S^2$. On the target side,  we will always think of $S^2$ as the unit sphere in $\mathbb{R}^3$, and consequently $S^2\times S^2 \subset \mathbb{R}^3 \times \mathbb{R}^3$.
We always work with 64-bit floating point numbers.

We now describe how to parametrize such functions, the main distinction is the depth, i.e. how many layers of neural networks we use to describe the function. 
We begin with simplest depth 0 case.

\subsection{Depth 0: only Fourier modes}

 Hence we will view $u$ as a map with domain $S^2$ or $T^2$, and target $\mathbb{R}^3$ but with image forced to lie on $S^2\subset \mathbb{R}^3$ (or $\mathbb{R}^6$ and image in $S^2\times S^2$). 

We will represent \[u=\mathcal{N}\circ \mathcal A\circ \mathcal{F}.\] Here $\mathcal{F}$ is a Fourier layer that maps from $X\rightarrow \mathbb{R}^{d}$, where $d$ is the number of Fourier modes we are considering.
In the case where $X=T^2$ with coordinates $(x,y)$, we can represent said Fourier expansion as
\[
\mathcal{F}(x,y) = \begin{pmatrix}
\sin x\\
\cos x\\
\sin y\\
\cos y\\
\sin 2x\\
\vdots\\
\sin x \cos y\\
\vdots
\end{pmatrix}
\]

In the case where $X=S^2$, we take the analogous expansion using spherical harmonics\footnote{In producing the first five leaves in Section~\ref{ss:gromovFoliation}, the constant spherical harmonic $Y_0^0$ was also included; this is redundant with the bias term and does not change the final function.}
\[
\mathcal{F}(x^1,x^2,x^3) = \begin{pmatrix}
Y_{-1}^{1}(x^1,x^2,x^3)\\
Y_{0}^{1}(x^1,x^2,x^3)\\
Y_{1}^{1}(x^1,x^2,x^3)\\
Y_{-2}^{1}(x^1,x^2,x^3)\\
\vdots
\end{pmatrix}.
\]
Here $x^1,x^2,x^3$ are the coordinates on $\mathbb{R}^3$ which contains $S^2$ as the unit sphere.
Here $Y_i^j$ are functions with input $\mathbb{R}^3$ and output $\mathbb{R}$ called the solid harmonic functions. 
They recover the usual spherical harmonic functions once we restrict to $S^2$. 
We use solid harmonic function for later ease of computation. 
When we train, plot and work with the function $u$, we only use the values of $Y_i^j$ on $S^2$.

The map $\mathcal A$ is an affine-linear map $\mathbb{R}^{d}\to \mathbb{R}^n$ where $n=3$ or $6$.
We represent each component of $\mathcal A$ as
\begin{equation}
    \mathcal A_{i}(z_j)=b_i+\sum_{j=1}^d w_{ij}z_j
\end{equation}
for to-be-determined real constants $w_{ij}$ (the \emph{weights}) and $b_i$ (the \emph{biases}).

Then the normalization layer $\mathcal{N}$ forces the function to lie on the target, either $S^2\subset \mathbb{R}^3$ or $S^2\times S^2 \subset \mathbb{R}^6$, so it is a function that takes the form
\begin{equation*}
    \mathcal N(z)=\frac{z}{\|z\|}\text{ or }\left(\frac{z_1}{\|z_1\|},\frac{z_2}{\|z_2\|}\right).
\end{equation*}

\subsection{Higher depth neural networks}

Because a general almost complex structure $J$ gives rise to a quasilinear elliptic PDE, it is natural to ask if including more nonlinearity in the representation of the function could improve the performance.
Use of deep neural networks to solve PDEs has been done in~\cite{raissi2019PINNs,deepritz}. 
The basic premise of a neural network is to write a potentially complicated function in terms of compositions of many much simpler nonlinear functions.
For definiteness we give the explicit formula for a depth 2 neural network. 
The construction of a depth 1 neural network is analogous.

As before, we take a Fourier feature layer $\mathcal{F}$ using Fourier modes (either sines and cosines or spherical harmonics depending on the source manifold).
In addition to providing a point of comparison to the natural model of the previous subsection, Fourier features provide added benefits in deep learning.
Fourier features in the initial layer of a neural network were popularized in~\cite{tancik2020fourierFeatures} in the context of mitigating spectral bias, the tendency of neural networks to preferentially learn low-frequency features  over high-frequency features~\cite{rahaman2019spectralBias}.
It serves a dual purpose for us, also eliminating the need to enforce periodic boundary conditions as a soft constraint inside the loss function.
Increasing the bandwidth of the Fourier feature layer allows resolving sharper features, although in practice this also enables the appearance of sharp singularities that are undesirable, see Table~\ref{tab:hyperparametersweep}.

We pick hidden layer widths $d_2\gg 1$ and $d_3\gg 1$; $d_4 =3$ or $d_4=6$ depending on whether the target is $S^2$ or $S^2\times S^2$; let $d_1$ denote the number of Fourier modes.

Our neural network then writes $u$ to be of the form
\begin{equation}
    u=\mathcal N\circ  \mathcal A^3\circ \sigma\circ \mathcal A^2\circ \sigma \circ \mathcal A^1\circ  \mathcal F,
\end{equation}
where
\begin{equation*}
    \mathcal A^i:\mathbb R^{d_i}\to \mathbb{R}^{d_{i+1}}
\end{equation*}
is an affine-linear map 
\begin{equation*}
    z^{i+1}_j =\mathcal A^i_j(z^i) = \sum_{k=1}^{d_i} W_{jk}^i z^i_k + b_j^i
\end{equation*}
parameterized by weights $W_{jk}^i$ and biases $b_j^i$ that will be learned during training.

A nonlinear activation function, the Gaussian-error linear unit (GELU)~\cite{hendrycks2016gelu}, approximated for ease of computation as
\begin{equation*}
    \sigma(z)=0.5z\left(1+\tanh\left[\sqrt{\frac{2}{\pi}}\left(z+0.0044715z^3\right)\right]\right),
\end{equation*}
is applied to the outputs of $\mathcal A^{1}$ and $\mathcal A^{2}$.
This function is a smooth, nonzero-curvature version of the ReLU activation function
\begin{equation*}
    \operatorname{ReLU}(z)=\max(0,z).
\end{equation*}
As before
\begin{equation*}
    \mathcal N(z)=\frac{z}{\|z\|}\text{ or }\left(\frac{z_1}{\|z_1\|},\frac{z_1}{\|z_2\|}\right)
\end{equation*}
is a normalization layer to ensure that the output lies on the target manifold (either $S^2 \subset \mathbb{R}^3$ or $S^2\times S^2 \subset \mathbb{R}^3\times \mathbb{R}^3$).

\section{The loss function}

Now that we have successfully parametrized a smooth function, we develop the tools to make it $J$-holomorphic.
To enforce $J$-holomorphicity, we have the standard loss for the relevant $\overline\partial $ operator. We fix metrics on both the domain and the target manifolds, and take
\begin{equation}
    \mathcal L_{\mathrm{cr}}=\int_{\Sigma}\ \|J(u)\circ du-du\circ j\|^2 d\textup{vol}\label{e:crloss}.
\end{equation}

The metric we take on the domain will always be the standard metric, whereas the metric on the target depends on $J$ as

\begin{equation}
    g_J(u,v)=\frac{1}{2}\left(\omega(u,Jv)-\omega(Ju,v)\right).
\end{equation}

We note $J$-holomorphic maps often come in high dimensional families. Assuming the almost complex structure $J$ is regular\footnote{Regular here is a technical term which means all moduli spaces of $J$-holomorphic curves are transversely cut out. This is called transversality in symplectic geometry, and in general establishing transversality is a delicate business. However, in all of the cases we consider the setup is simple enough all of the transversality conditions we need are satisfied, see for instance \cite{wendl}.}, for a fixed (closed) Riemann surface $(\Sigma,j)$ of genus $g$, the dimension of $J$-holomoprhic maps from $\Sigma \rightarrow M$ is a space of dimension
\[
n(2-2g) +2c_1(u^*TM)
\]
by the Riemann-Roch formula. 
We fix a collection of points $\{x_1^*,..,x_k^*\}$ points on the domain $(\Sigma,j)$, and a collection of disjoint submanifolds $\{N_1,...N_k\}$ on $M$ (for us we will mostly these to be points) on the target symplectic manifold to impose the constraint $u(x_i^*)\in N_i$. The curves that satisfy these constraints will live in a space of dimension (assuming transversality holds)
\begin{equation}
n(2-2g) +2c_1(u^*TM) -\sum_i\textup{codim}(N_i).\label{e:markedpointdimformula}
\end{equation}
We usually choose the marked point constraints so that the above dimension is zero, so that given a $J$-holomorphic curve satisfying the above constraint, there are no nearby curves satisfying the same constraint.
\begin{remark}\label{rmk:holmaps}
Usually to obtain meaningful counts of $J$-holomorphic curves and obtain the required transversality conditions, we need to let the domain complex structure $j$ vary if we have a curve of genus $\geq 1$, and quotient by the automorphism groups of the domain. This sometimes results in a smooth orbifold whose dimension is slightly different from the formula we gave above \cite{wendl}. However in the cases we consider the set up is simple enough the needed transversality  is automatic \cite{wendl} even with fixed domain complex structure, and for a better fit with our numerical algorithms we do not quotient by automorphisms of the domain. To be completely precise about terminology, technically speaking we are considering \emph{parameterized} $J$-holomorphic curves instead of $J$-holomorphic curves. However we gloss over this point via abuse of notation.
\end{remark}

To numerically enforce the marked point constraints, we add a marked point component to the loss. Let $\{x_i^*\}$ denote a fixed collection of marked points on the domain $(\Sigma,j)$ and let $N_i$ be a collection of disjoint submanifolds on the target $M$ ,
\begin{equation}
\mathcal{L}_{\mathrm{mp}}=\sum_{i=1}^{n_{\mathrm{mp}}} \operatorname{dist}\left(u(x^*_i), N_i\right)^2.\label{e:ancloss}
\end{equation}

Primarily for diagnostic purposes (and to add a slight guardrail during the beginning phase of training when the learning rate is high), when the target is $S^2$ we consider the degree loss of a map $u:X\rightarrow Y$ with topological degree $d$ as 
\begin{equation}
    \mathcal L_{deg}=\left(\int_{X} u^*\omega-d\right)^2
    \label{e:degloss}
\end{equation}
where $\omega$ is the symplectic form on the target. If the target is $S^2\times S^2$ with the symplectic form $\omega_1+\omega_2$, we consider the integral of $\omega_1$ and $\omega_2$ separately as the bidegree. This degree should be invariant through the training process as it is topological. It will be important for us later when we diagnose during if a component of the curve has bubbled off.

To minimize the appearance of concentrated spikes in energy (which is a known degeneration mechanism for $J$-holomorphic curves), we also add a regularization term. Let $\{x_i\}$ denote a random collection of sample points on the domain,
\begin{equation}
    \mathcal{L}_{\mathrm{reg}}=\sum_{i=1}^{N} \left(|\nabla u(x_i)|-\tau\right)_+^2,\label{e:gradregularization}
\end{equation}
where $v_+$ denotes the positive part of $v$, penalizing points where the gradient of $u$ is larger than the threshold $\tau$.
Instead of enforcing a universal upper limit on the derivative (which we cannot know a priori), we compute $\tau$ adaptively based on the ninety-ninth percentile of pointwise gradients of $u$ currently during the training process; this means this term will be positive, but we can schedule the weight to turn off near the end of training.

The final loss function will be a weighted linear combination of the above individual functions
\begin{equation}
    \mathcal L[u]=w_{\mathrm{cr}}\mathcal L_{\mathrm{cr}}+w_{\mathrm{deg}}\mathcal L_{\mathrm{deg}}+w_{\mathrm{mp}}\mathcal L_{\mathrm{mp}}+w_{\mathrm{reg}}\mathcal L_{\mathrm{reg}},
\end{equation}

To approximate the integrals here numerically, we sample a large number of points on the domain and estimate the integral numerically by taking a large batch size (a few thousand to over ten thousand points) in a grid adapted to the source manifold.
Randomness in this sampling can mitigate overfitting and assist in the training process.

The choice of weights 
$w_{\mathrm{cr}}$, $w_{\mathrm{deg}}$, $w_{\mathrm{mp}}$, and $w_{\mathrm{reg}}$ are chosen to vary during the training phase. How these weights vary will depend on the precise training schedule we use, and as we shall see, they will make noticeable differences in training dynamics.

Another parameter we will track but not use for training is what we call the bubbling ratio. 
In cases where the holomorphic curve we are trying to find is explicitly known, we wish to analyze the trials that did not converge to the known curve and understand why they did not converge. 
As we shall see some of them failed to converge due to local high concentrations of energy. 
Let $u_{\mathrm{GT}}$ be the known holomorphic function we are trying to find, and let $u_1$ denote the non-holomorphic function our methods produced, then we define
\[
\textup{bubbling ratio} = \frac{\|du_1\|^2_{L^\infty}}{\| du_{\mathrm{GT}}\|^2_{L^\infty}}.
\]
In the above, we compute $\|du_1\|^2_{L^\infty}$ by choosing a coordinate chart in the domain and target, and compute $L^\infty$ of the Hilbert-Schmidt norm of the matrix $du_1$. 

\section{The training paradigms}\label{sec:training}

We now describe the main paradigms we shall use for our training. The input will always be a smooth function, and desired output will be a $J$-holomorphic curve satisfying the marked point constraints.

\subsection{Paradigm 1}

\begin{figure}[H]
\begin{center}
\includegraphics[width=0.5\linewidth]{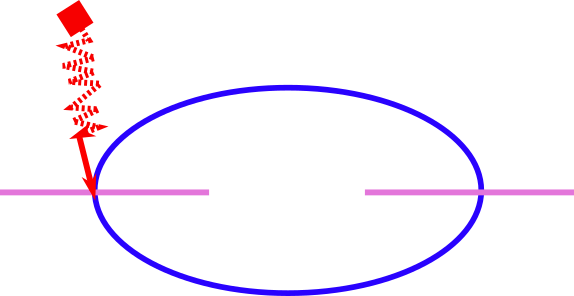}
\end{center}
\caption{The blue locus is the space of $J$-holomorphic curves. The magenta locus is the space of smooth maps that satisfy the marked points constraints. The $J$-holomorphic curve we want lies at the intersection of the blue and magenta region. The red square denotes the initial smooth curve we started with, and training in this paradigm means we set up the loss function that steers it directly towards the intersection, taking account into the Cauchy-Riemann loss and the marked point loss simultaneously (along with the regularization term). The jagged dash line indicates the noisy Adam gradient descend phase, and the filled straight line indicates the polish phase.}\label{f:paradigm1}
\end{figure}

This approach is a straightforward schedule for finding a single curve; a cartoon illustration is given in Figure~\ref{f:paradigm1}.

\begin{enumerate}
    \item The first step is pretraining to fit a smooth function $\tilde u_0$ in the right homology class; this is skipped in certain steps in Sections~\ref{ss:gromovFoliation} and~\ref{ss:bifurcation} when we are starting from a curve that is already close to the desired one. 

    The pretrained model $u_0$ is a least-squares regression to $\tilde u_0$.
    In the depth-zero model, this is done by linear regression.
    In higher-depth PINNs, this is done using the Adam optimizer~\cite{diederik2015adam}, a first-order gradient descent method. This a stochastic gradient descent optimizer built by the machine learning community 2015, with features like adaptive moment estimation designed to search for the solution even when the loss landscape is not strictly convex. One of its advantages is robustness to large parameter counts in terms of cost.

    We split the main training into two phases, a globalizing Adam phase and a polishing quasi-Newton phase. In both phases we turn on both $\mathcal{L}_{\mathrm{cr}}$ and $\mathcal{L}_{\mathrm{mp}}$ while adjusting their respective weights (and optionally $\mathcal{L}_{\mathrm{reg}}$). 
    We keep the learning rate small enough so that homology class does not jump during training, and this is monitored by $\mathcal{L}_{\mathrm{deg}}$ (even if $w_{\mathrm{deg}}\equiv 0$, the degree integral is computed to confirm topology has not degenerated).

    \item In the Adam phase, the sample points in computing $\mathcal{L}_{\mathrm{cr}}$ are randomly redrawn each iteration, providing additional smoothing.
    The weights $w_{\mathrm{deg}}$ and $w_{\mathrm{mp}}$ start large and are scheduled to decay as $\mathcal{L}_{\mathrm{cr}}$ hits predefined thresholds.
    The gradient regularization weight $w_{\mathrm{reg}}$, if present, is ramped up only in the middle of this phase to give the network time to approach the correct solution and is annealed off by the end to ensure that the final formula for $\mathcal L$ in this phase gives $0$ for a $J$-holomorphic curve.
    The goal of the weight schedules is to keep a balance where the weighted Cauchy-Riemann loss remains around ten to a hundred times larger than the weighted marked-point loss.
    Learning rate starts at $5\cdot 10^{-4}$.
    After a set number of iterations where $\mathcal L_{\mathrm{cr}}$ does not improve from its best value so far, the learning rate is halved, with $5\cdot 10^{-7}$ as its floor.

    \item In the polish phase, the weights are optionally rebalanced. The sample grids for integrals are fixed, and a second-order quasi-Newton method is used.
    In the depth-two case, we use L-BFGS~\cite{liu1989lbfgs}, a quasi-Newton method that estimates required second derivative information using the recent history and is optimized for limited memory consumption even with large numbers of parameters.
    We use a fixed grid in the source manifold that is chosen densely enough to avoid aliasing artifacts.
\end{enumerate}

\subsection{Paradigm 2}

\begin{figure}[H]
\begin{center}
\includegraphics[width=0.5\linewidth]{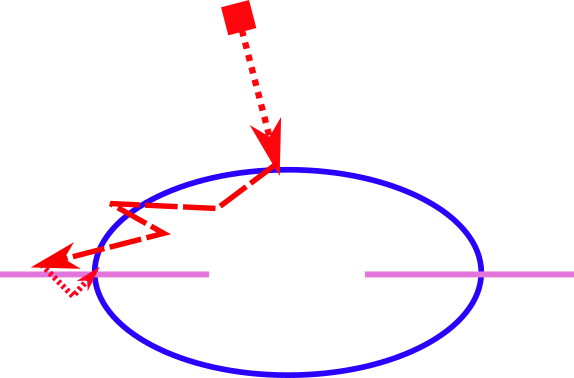}
\end{center}
\caption{The blue locus is the space of $J$-holomorphic curves. The magenta locus is the space of smooth maps that satisfy the marked points constraints. The $J$-holomorphic curve we want lies at the intersection of the blue and magenta region. The red square is the beginning of the smooth curve. We use different lines to indicate different phases of training. We first descend as quickly as possible to being $J$-holomorphic without any regard for marked point constraints. Then we drag the curve to the magenta region along the blue locus by ramping up the marked point loss and using the Adam optimizer, fitting the marked point constraints better at the cost of increased Cauchy-Riemann loss. When it gets close enough we polish as a third phase to bring down the overall loss.}\label{f:paradigm2}
\end{figure}

In practice, the Cauchy-Riemann loss $\mathcal{L}_{\mathrm{cr}}$ and the marked point loss $\mathcal{L}_{\mathrm{mp}}$ can be difficult to minimize simultaneously. See the discussion in Section \ref{sss:hyperparametersweep} for examples where a fast-training depth-0 model fails to converge. 
To alleviate the difficulties associated with the highly nonconvex loss landscape, paradigm 2 partially separates the two concerns.
A cartoon drawing is given in Figure~\ref{f:paradigm2}.

We use this paradigm only with depth-0 models because of the cost associated with quasi-Newton optimization when parameter counts are large.

\begin{enumerate}
    \item The first step, as in paradigm 1, is pretraining to fit a smooth function $\tilde u_0$ in the right homology class. 

    \item Then, we use the Levenberg-Marquardt algorithm \cite{Levenberg,Marquardt}, a quasi-Newton method that becomes more like gradient descent when steps are not as successful, to reduce $\mathcal{L}_{\mathrm{cr}}$.
    We  use $w_{\mathrm{mp}}=0$, $w_{\mathrm{cr}}=100$, and $w_{\mathrm{deg}}=10$, focusing purely on making the curve $J$-holomorphic.
    The resulting curve will generally be far from satisfying the marked point constraints, but it will be very close to the manifold of $J$-holomorphic curves.

    \item To move along the manifold of 
    $J$-holomorphic curves towards the locus where the marked points constraints are satisfied, we use Adam with fixed, positive weights on both the marked point and Cauchy-Riemann loss.
    We ramp $w_{\mathrm{mp}}$ up from $0$ to $1000$ (with $w_{\mathrm{cr}}=100$) and leave $w_{\mathrm{deg}}=0$.
    This will result in several orders of magnitude worse $\mathcal{L}_{\mathrm{cr}}$, wandering slightly away space of $J$-holomorphic curves, but in exchange we greatly reduce the marked point loss, bringing the curve much closer to satisfying the point constraints. 

    \item Once we are near the desired solution, a final quasi-Newton polish phase (with the same weights as at the end of the Adam phase) returns us to the manifold of $J$-holomorphic curves and tightens the fit of the marked point constraints.
\end{enumerate}

\section{Validation with integrable $J$}

\subsection{The Weierstrass function on the torus} \label{sec:wp}

Consider the torus $T^2=S^1\times S^1$ with the product almost complex structure $j$ (i.e. it interchanges the two $S^1$ factors). There is a well known meromorphic function called the Weierstrass $\wp$-function which we can view as a degree 2 map to the Riemann sphere, $\wp: T^2\rightarrow S^2$. 

If we view $T^2=\mathbb{R}^2/\Lambda$, where $\Lambda=\mathbb{Z}^2$, then the Weierstrass $\wp$ function is given by
\[
\wp(z) = \frac{1}{z^2} +\sum_{\lambda \in \Lambda\setminus 0}\left ( \frac{1}{(z-\lambda)^2}-\frac{1}{\lambda^2}\right )
\]

As our first proof of concept we use our numerical algorithm to find this function. We start with an arbitrary degree 2 smooth function $f:T^2\rightarrow S^2$, impose 4 generic marked point constraints, and apply our methods to recover a holomorphic map satisfying the marked point constraints.

Here we are exploring the space of parametrized degree 2 holomorphic maps $T^2\rightarrow S^2$. This space is explicitly known. It is 8 dimensional, and all elements of this space are given by $f(z) = h\circ \wp \circ k(z)$, where $k(z)$ is a translation on $T^2$ and $h(z)$ is a fractional linear transformation.

When we impose the 4 marked point constraints on this space of holomorphic maps, we find for generic choice of marked points there are actually \emph{two} different maps from $T^2\rightarrow S^2$ satisfying those marked point constraints. We will see both of these functions make an appearance in the tests below.

For ease of visualization since it moves the interesting features away from the edges of our diagrams, we work with the translate
\begin{equation*}
    \wp^{(1)}(z)=\wp(z+0.31+0.42i).
\end{equation*}
We choose the marked point constraints such that $u$ must take the same values as $\wp^{(1)}$ at the four points
\begin{equation*}
    0.2+0.3i,\, 0.6+0.2i,\, 0.35+0.7i,\, \text{and } 0.8+0.85i.
\end{equation*}
This uniquely determines the second holomorphic curve, which we denote by $\wp^{(2)}$.
The translation and fractional linear transformation taking $\wp^{(1)}$ to $\wp^{(2)}$ are given approximately by
\begin{align*}
    k(z)&=z+0.905+ 0.135i\\
    h(z)&=\frac{(1.648-4.156i)z +( -12.126 + 28.503i)}{(0.365 + 0.673i)z+1}
\end{align*}
as found by a Newton solve.

\subsubsection{Results with Paradigm 1}\label{sss:paradigm1}

We use a depth 2 PINN with Fourier radius $8$ to parameterize our functions (as suggested by the hyperparameter sweep of Section~\ref{sss:hyperparametersweep}), and apply paradigm 1.

We start with 16 different smooth degree-$2$ functions in various configurations.
We found depending on the initial configuration, 6 of them converged to $\wp^{(1)}$ and 4 converged to $\wp^{(2)}$. The remaining 6 did not converge at all, and the results are summarized in Table~\ref{tab:wpparadigm1}. 
We also found that those initial smooth configurations that started out closer to satisfying the marked point constraints had better chances of convergence than those initial configurations that started out with greater marked point loss. 

\begin{remark}
    The fact we were able to converge to both $\wp^{(1)}$ and $\wp^{(2)}$ given different initializations opens up the possibility of using our methods with a large number of random initializations to enumerate (numerical) $J$-holomorphic curves in a symplectic manifold. A potential application of this would be numerical computations of Gromov-Witten invariants, which is of great interest to symplectic geometers.
\end{remark}

\begin{table}[H]
    \centering
    \begin{tabular}{c|c|c|c|c|c|c}
         final basin & \# & final $\mathcal{L}_{\mathrm{cr}}$ & max final mp dist &  $\operatorname{dist}(u,\wp^{(1)})$ &  $\operatorname{dist}(u,\wp^{(2)})$ & degree \\
         \hline 
         $\wp^{(1)}$ & 6 & $7\cdot 10^{-6}$ & $4.57\cdot 10^{-5}$ & $4.37\cdot 10^{-4}$ & $0.645$ & $2.08$ \\
         $\wp^{(2)}$ & 4 & $1\cdot 10^{-5}$ & $2.13\cdot 10^{-4}$ & $ 0.645$ & $1.3\cdot 10^{-3}$ & $1.92$\\
         $\mathcal{L}_{\mathrm{cr}}$ failures & 5 & $\approx 10^{-4}$ to $\approx 10^{-2}$&  $\approx 10^{-4}$ to $\approx 10^{-1}$ & all $>0.4$ & all $>0.4$  & $2.01$ \\
         deg failures & 1 & $6.21\cdot 10^{-2}$ & $7.53\cdot 10^{-3}$ & $0.461$ & $0.319$ & $0.764$*
    \end{tabular}
    \caption{Summary statistics from a multi-initialization run with $16$ randomly-oriented tubes using paradigm 1, averaged by which curve the run converged to or what kind of failure was observed.
    Note the degrees are slightly different from integer values since the degree integral is computed numerically.
    The relatively high variance between runs leads us to report fewer significant digits.}
    \label{tab:wpparadigm1}
\end{table}

We observe two distinct failure modes in the 6 runs that did not converge. 
The first class we call analytic failures where the $CR$ loss was still high after the Adam phase and did not go down (and sometimes even increased) in the polish phase. 
This suggests after the Adam phase the function landed far away from the actual solution, and L-BFGS was unable to help convergence because of this bad starting point. 
We note that each of these had the correct degree. 
The other class of failures we call topological failures where even the degree is not preserved. We suspect this is some form of bubbling phenomena analogous to the degeneration of $J$-holomorphic curves. 
The one case where topological failure occurred had a bubbling ratio estimated at $114$, which indicates great concentration of energy.

While it is an encouraging result that many initializations converged and we were able to find both $\wp^{(1)}$ and $\wp^{(2)}$, it is in our opinion very much worth exploring how to increase the convergence fraction, run the code faster (PINNs are expensive to run computationally), and reduce the final error for the converged runs. 
With that in mind we turn to paradigm 2 in the next subsection where indeed in some regimes we do see better convergence, much smaller Cauchy-Riemann loss, and an order of magnitude speed up. 
See also discussion in Section~\ref{sss:hyperparametersweep}, in particular the wall times in Table~\ref{tab:hyperparametersweep}.

\subsubsection{Results with Paradigm 2}\label{sss:paradigm2}
We use depth-0 and Fourier radius $8$ coupled with paradigm 2, keeping as many other aspects of the training as possible the same to isolate the differences arising from the depth and paradigm shifts. 
We start with the same 16 random smooth degree-$2$ functions as in Section~\ref{sss:paradigm1} above. 
We found after the phase 2 (the second order quasi-Newton method minimizing $\mathcal{L}_{cr}$ without any constraint in $\mathcal{L}_{pm}$), \emph{all} of the initial smooth guesses converged to holomorphic curves. We record that all 16 runs had Cauchy-Riemann loss below $6\times 10^{-9}$, which is extremely good convergence compared to the results of paradigm 1. Based on this comparison we hypothesize that the simple act of turning on the marked point loss greatly complicates the loss landscape in a way the prevents convergence (perhaps by losing convexity of the loss functional). We also remark in both paradigm 1 and 2 we had unconverged functions where the marked points constraints were actually satisfied all very well, but the Cauchy-Riemann loss refused to go down.

\begin{figure}[H]
    \centering
    \includegraphics[width=0.8\linewidth]{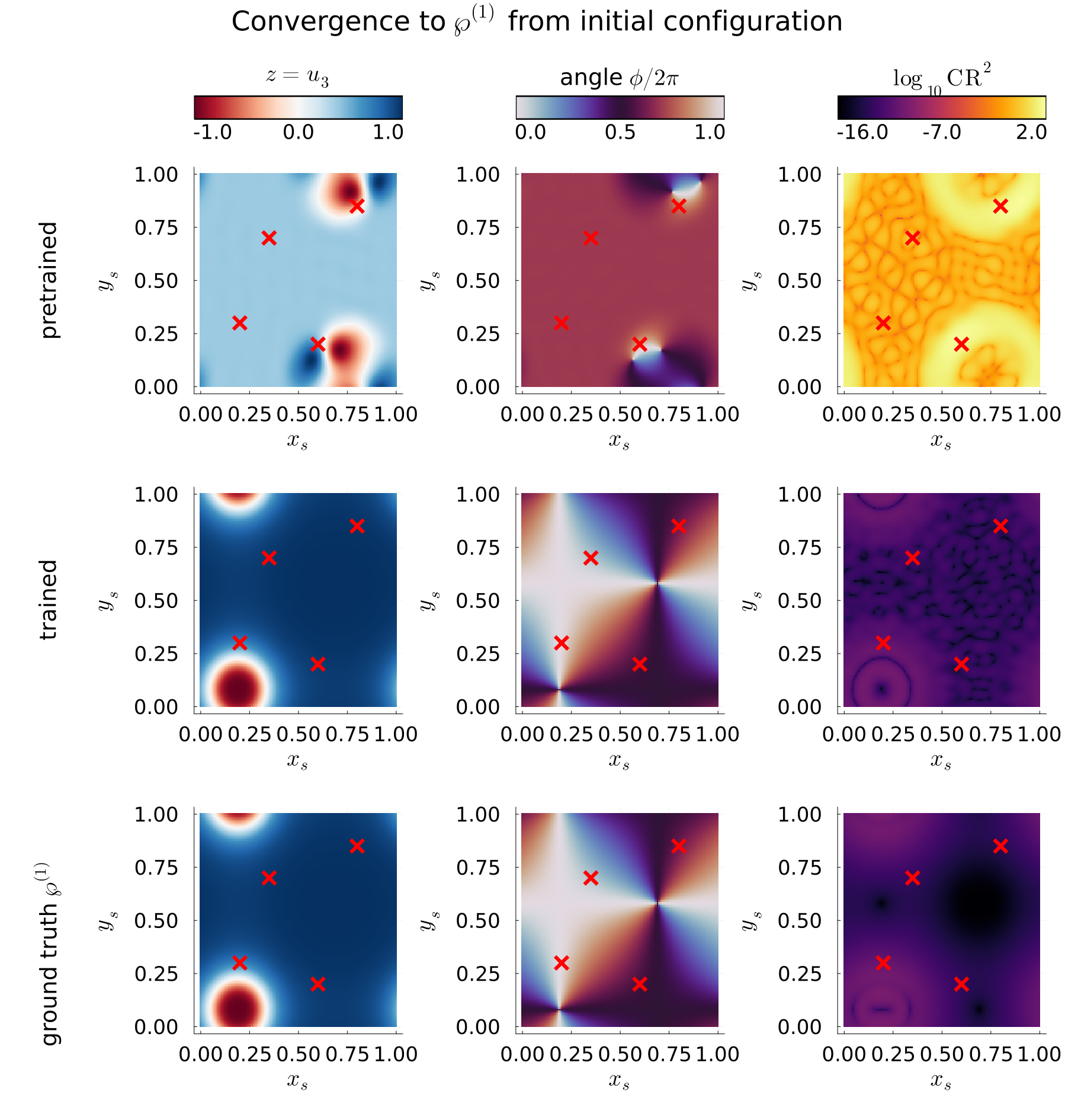}
    \caption{A representative initial pretrain target that converged to $\wp^{(1)}$.
    Each marked point is denoted with a red `x.'
    The subscript $s$ denotes source coordinate.}
    \label{fig:p2_wp1}
\end{figure}

\begin{figure}[H]
    \centering
    \includegraphics[width=0.8\linewidth]{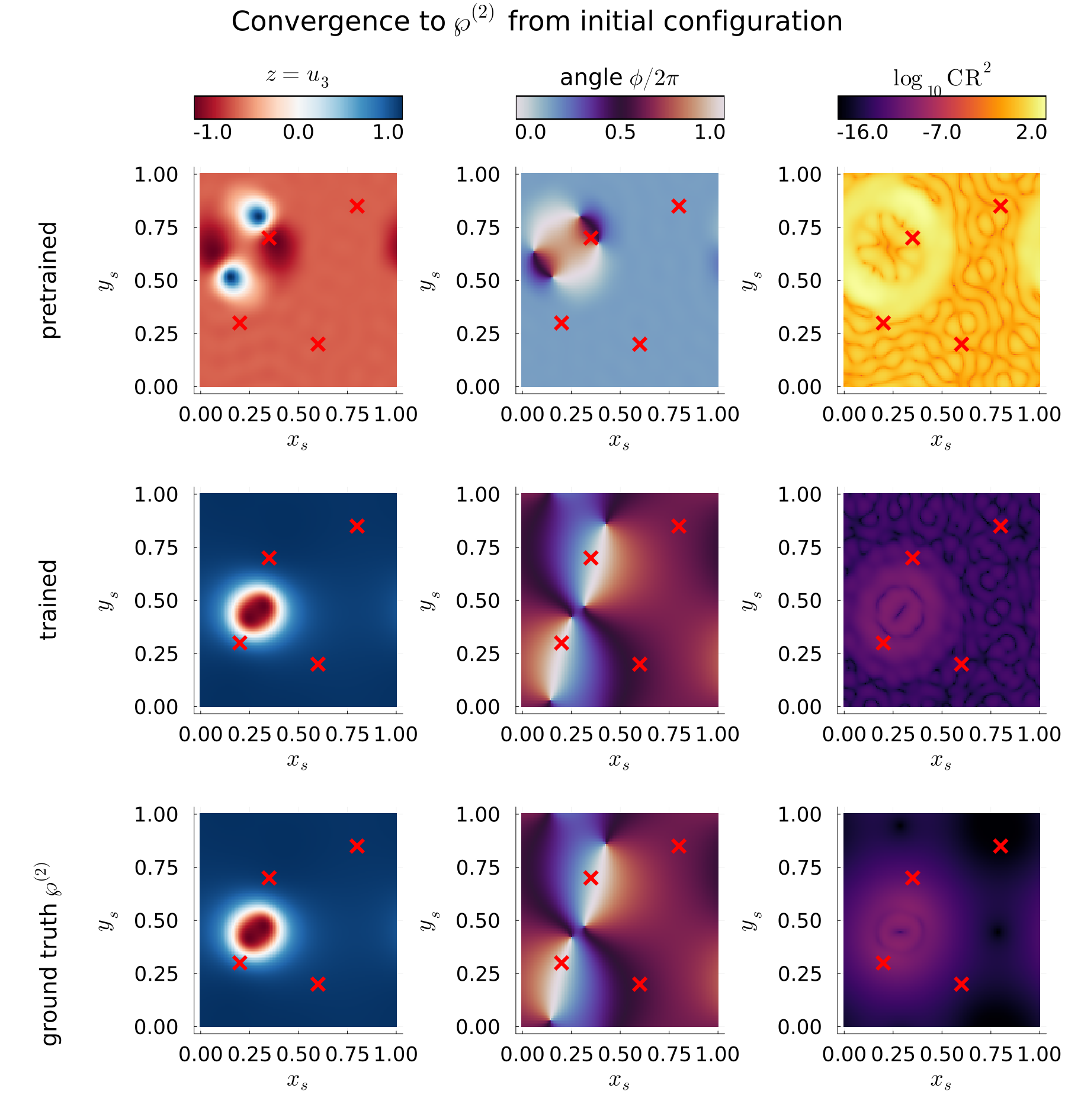}
    \caption{A representative initial pretrain target that converged to $\wp^{(2)}$}
    \label{fig:p2_wp2}
\end{figure}

As we gradually turned on the marked point constraints in phase 3, 11/16 holomorphic curves converged to $\wp^{(1)}$ or $\wp^{(2)}$. 
The data is summarized in table~\ref{tab:wpparadigm2}.  Plots of two successful runs are shown in Figures~\ref{fig:p2_wp1} and~\ref{fig:p2_wp2}.

We are visualizing the domain as $T^2=[0,1]^2$ which are the squares shown in the figures. 
We think of $S^2$ as the unit sphere in $\mathbb{R}^3$.
To encode a map from $g: T^2\rightarrow S^2$ we use the following convention. 
The first column is the $z$ component of $g$, and the second column is the angle of $g$ when projected to the $x-y$ plane. 
In other words we are using cylindrical coordinates to describe the target. 
The third column shows the (log of)  Cauchy Riemann loss, which tells us how $J$-holomorphic the function is.

The first row shows the smooth map $f$ we parametrized after pretraining.
$f$ is a degree 2 smooth map from $T^2\rightarrow S^2$ we wrote down by hand. The 2nd row shows the $J$-holomorphic curve after paradigm 2.  The 3rd row shows the actual plot of the known Weierstrass $\wp$ function suitably transformed to satisfy the marked points.

\begin{table}[H]
    \centering
    \begin{tabular}{c|c|c|c|c|c|c}
         final basin & \# & final $\mathcal{L}_{\mathrm{cr}}$ & max final mp dist &  $\operatorname{dist}(u,\wp^{(1)})$ &  $\operatorname{dist}(u,\wp^{(2)})$ & degree \\
         \hline 
         $\wp^{(1)}$ & 6 &$3.026 \times 10^{-12}$ & $1.94 \times 10^{-12}$ & $1.351 \times 10^{-7}$ & $0.6454$ & $2.077$\\
         $\wp^{(2)}$ & 5 & $7.72 \times 10^{-12}$ & $1.245 \times 10^{-11}$ & $0.6454$ & $1.214 \times 10^{-7}$ & $1.924$\\
         $\mathcal{L}_{\mathrm{cr}}$ failures & 3 & $9.196 \cdot 10^{-3}$ & $5.932 \cdot 10^{-2}$ &  -- & -- &  $1.955$ \\
         deg failures & 2 & -- & $5.875\cdot 10^{-8}$ & -- & -- & -- %58.556
    \end{tabular}
    \caption{Summary statistics from a multi-initialization run with $16$ randomly-oriented tubes using paradigm 2, broken down by which curve (if any) the run converged to. In this table where we report averages it indicates all runs had comparable values.
    More significant figures are given as compared to Table~\ref{tab:wpparadigm1} because the final learned curves clustered together more closely with paradigm 2.}
    \label{tab:wpparadigm2}
\end{table}
We observe because we used depth 0 networks instead of depth 2, the runs that converged ran faster by an order of magnitude.
Having access to a more powerful (but slower per-parameter) optimizer helped ensure that for the converged runs, the Cauchy-Riemann loss was smaller by 6 orders of magnitude.

We again notice the two distinct failure modes, analytical and topological. All cases in the analytical failure mode had the right degree.
The two cases in topological failure modes both had reported degrees approximately 3. 
All unconverged runs has a bubbling ratio above $2$, and for the runs classified as ``topology failures,'' the bubbling ratios were exceptionally bad, estimated at $8.42$ and $6260$.

Note if we take all the results of paradigm 1 and paradigm 2 combined, 14 out of 16 initializations converged. For seeds where both paradigms converged, they converged to the same function.

\subsubsection{Hyperparameter sweep} \label{sss:hyperparametersweep}

For the case of $\wp$ we compare in detail the pros and cons of multi-layer neural-networks vs depth 0 models. 
We ran the code with depth 0, 1, and 2 and varied both the width of the hidden layers and the number of Fourier modes, holding the schedule constant. 
We implemented paradigm 1 in which we include both $\mathcal{L}_{\mathrm{cr}}$ and $\mathcal{L}_{\mathrm{mp}}$. 
The sweep was run using an Apple M4 Pro with 24 GB of RAM.
We note that all experiments in this paper were done either with that 2024 MacBook Pro or with other (slower) consumer laptops.

We found for depth-0 networks, none of the 3 trials converged. And almost universally, convergence was worse the more Fourier modes we added. 
We hypothesize the fact including more Fourier modes allowed the depth 0 models to capture higher frequency phenomena more easily which led to the bubbling off of the curves. We found as we added more layers in the neural network the bubbling got better and we converged, a phenomena which warrants additional experiments and explanation.

% max $\mathcal{L}_{\mathrm{cr}}$
\begin{table}[H]
    \centering
    \begin{tabular}{c|c|c|c|c|c|c|c|c}
         depth & width & fourier radius & wall time (s) &  $\mathcal{L}_{\mathrm{cr}}$ & $\left\|\overline{\partial}u\right\|_{L^\infty}^{2}$ & degree & bubbling ratio & basin\\
         \hline
         $0$ & -- & $12$ & $4205.3$ & $0.14574$ & $22.388$ & $2.0001$ & $2.8976$ & ?\\
$1$ & $512$ & $12$ & $9837$ & $0.00015706$ & $0.36595$ & $2.0$ & $1.0024$ & $\wp^{(1)}$\\
$1$ & $1024$ & $12$ & $15172.0$ & $0.0005411$ & $2.6717$ & $2.0001$ & $1.09$ & $\wp^{(1)}$\\
$2$ & $512$ & $12$ & $15720.0$ & $0.00010756$ & $0.18312$ & $2.0$ & $1.0098$ & $\wp^{(1)}$\\
$2$ & $1024$ & $12$ & $29736.0$ & $7.4223 \times 10^{-5}$ & $0.12466$ & $2.0$ & $1.0069$ & $\wp^{(1)}$\\
$0$ & -- & $18$ & $6126.9^*$ & $0.1605$ & $18.295$ & $2.0001$ & $4.9704$ & ?\\
$1$ & $512$ & $18$ & $15446$ & $0.074002$ & $274.18$ & $2.0735$ & $57.794$ & ?\\
$1$ & $1024$ & $18$ & $23220.0$ & $0.014226$ & $31.668$ & $1.9967$ & $38.491$ & ?\\
$2$ & $512$ & $18$ & $21426.0$ & $0.0010867$ & $1.2829$ & $2.0$ & $2.9306$ & ?\\
$2$ & $1024$ & $18$ & $39685.0$ & $0.00040417$ & $0.26142$ & $2.0$ & $1.0983$ & $\wp^{(1)}$\\
$0$ & -- & $24$ & $6504.1^*$ & $0.16312$ & $34.98$ & $2.0007$ & $6.3515$ & ?\\
$1$ & $512$ & $24$ & $25701$ & $0.4305$ & $1579.0$ & $1.9903$ & $67.631$ & ? \\
$1$ & $1024$ & $24$ & $35069.0$ & $0.0040935$ & $1.5129$ & $1.9998$ & $4.9012$ & ?\\
$2$ & $512$ & $24$ & $29024$ & $0.0073194$ & $13.471$ & $2.0002$ & $2.8304$ & ?\\
$2$ & $1024$ & $24$ & $55141.0$ & $0.0075857$ & $11.499$ & $2.0$ & $3.439$ & ?
    \end{tabular}
        \caption{Results of the parameter sweep.
        Runs marked * terminated early in the polish phase due to non-improvement
        For depth-0 networks, there is no hidden layer for the width parameter to affect.}
    \label{tab:hyperparametersweep}
\end{table}

\begin{figure}[H]
    \centering
    \includegraphics[width=0.9\linewidth]{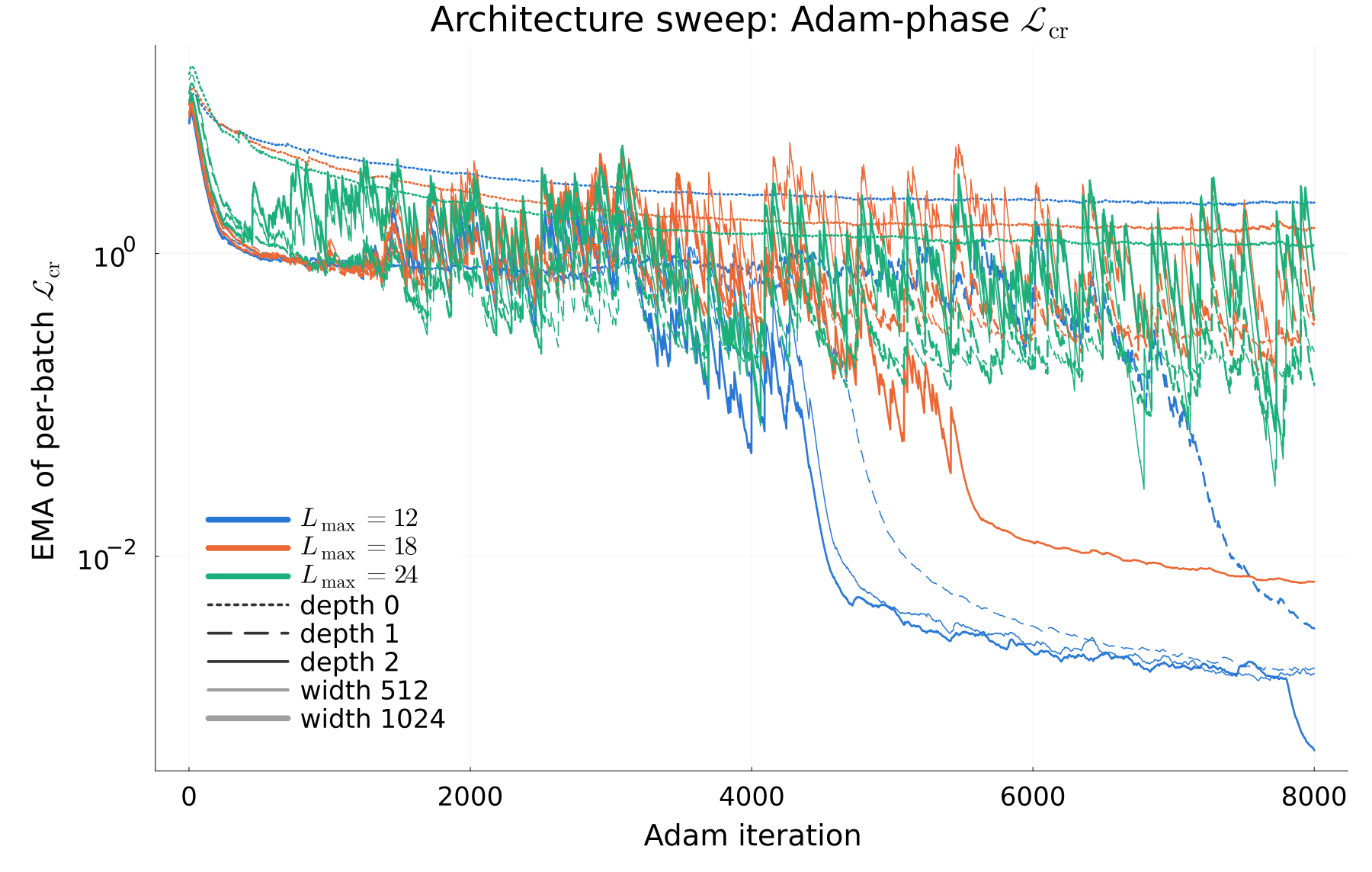}
    \caption{The loss curves for the Adam phase of the hyperparameter sweep.
    Note that most  have plateaued by the end of the alloted 8000 iterations and that the lowest loss belongs to depth $2$, width $1024$, Fourier radius $12$.
    This motivates our hyperparameter choices for paradigm 1.}
    \label{fig:hyperparamatersweep_losscurves}
\end{figure}

\subsection{A bidegree-$(5,3)$ holomorphic map $S^2\to S^2\times S^2$}\label{sec:5-3}
As our next proof of concept we consider the target to be $S^2\times S^2$ with the product complex structure $J_0$. We consider a $J$-holomorphic curve $u:S^2\rightarrow S^2\times S^2$, and the map is given by
$u(z)=\left (z^5,h(h^{-1}(z)^3)\right)$, where $h$ is fractional linear transformation coming from rotation in $\mathbb{R}^3$. 
We run paradigm 1 with depth 2 PINN.

We again start with a smooth function $f:S^2\rightarrow S^2\times S^2$ with the same homology class but intentionally far from the ground truth. 
To make sure we end at the same curve we place 18 marked points on the domain; 7 are used to fix the degree 3 component and 11 are used to fix the degree 5 component.
The marked points used to fix each component are distributed across the sphere to ensure the optimizer is well conditioned.
We place the marked points for each factor by starting from one randomly-chosen point, then placing points to maximize the minimum distance from the existing marked points, ensuring broad coverage of the domain.

We plot the graph of our output as follows. We put spherical coordinates on both the domain $S^2$ ($\phi,\theta$) and target $S^2\times S^2$ $(\phi_1,\theta_1,\phi_2,\theta_2)$. This means to specify a given map we need 4 images. We then plot the Cauchy-Riemann loss which measures how $J$-holomorphic the map is.

In the first row we plot the pretrained smooth function, and see it has large Cauchy Riemann loss. We next plot the curve after training in the second row, for which we see greatly reduced CR loss. We note there is some distortion of the plot near the top and bottom edges, but this is an artifact of us choosing spherical coordinates. Finally we plot the known holomorphic function in the last row, and see that the two graphs are close to each other.

\begin{figure}[H]
    \centering
    \includegraphics[width=0.9\linewidth]{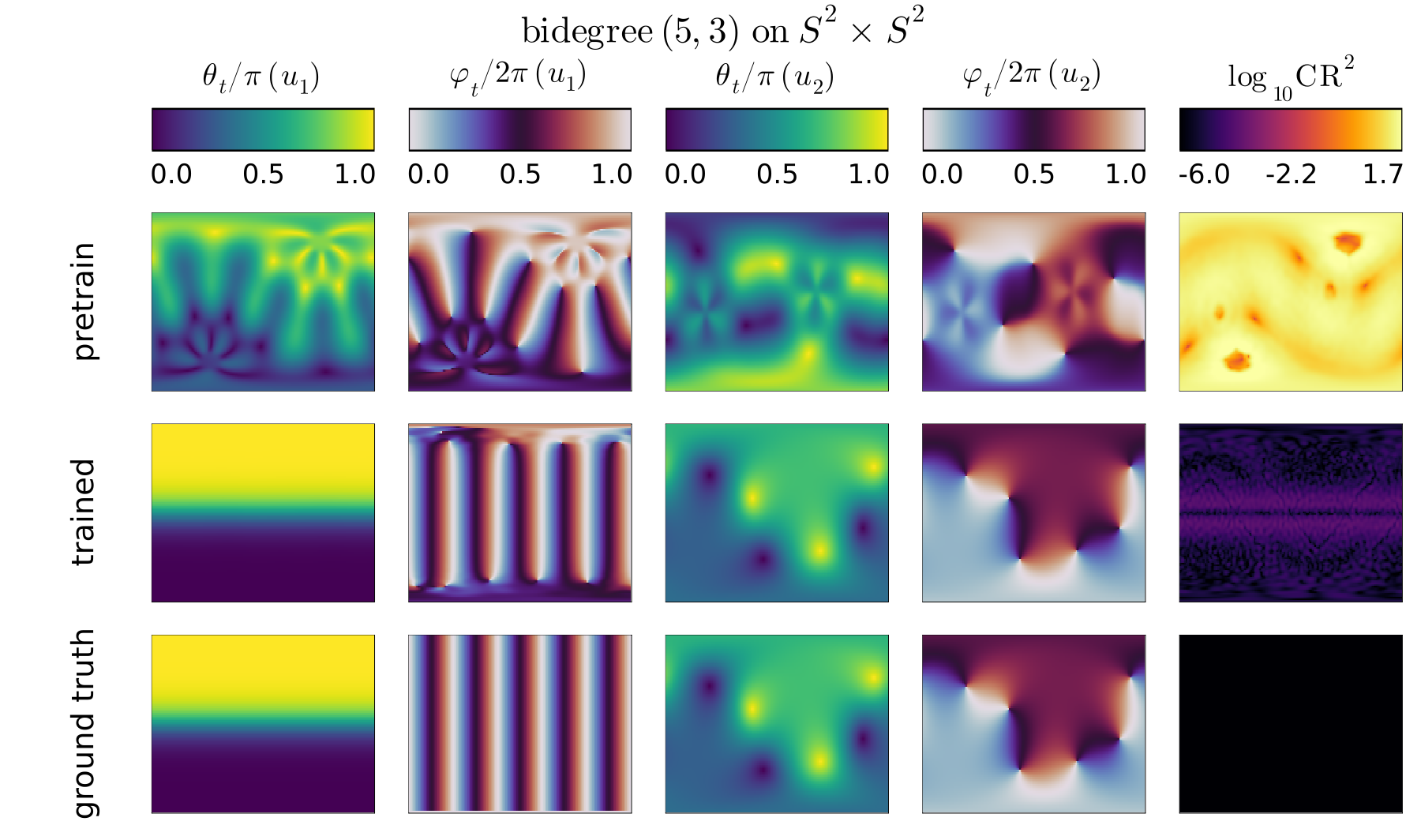}
    \caption{The image and $\mathcal{L}_{\mathrm{cr}}$ for the pretrained and trained networks as compared to the ground truth.
    Note the additional twisting in the smooth, non-holomorphic pretrain that is smoothed out out by the training.
    The subscript $t$ denotes a coordinate in the target manifold.
    }
    \label{fig:5_3_image}
\end{figure}\label{fig:5-3}

\begin{figure}[H]
    \centering
    \includegraphics[width=0.7\linewidth]{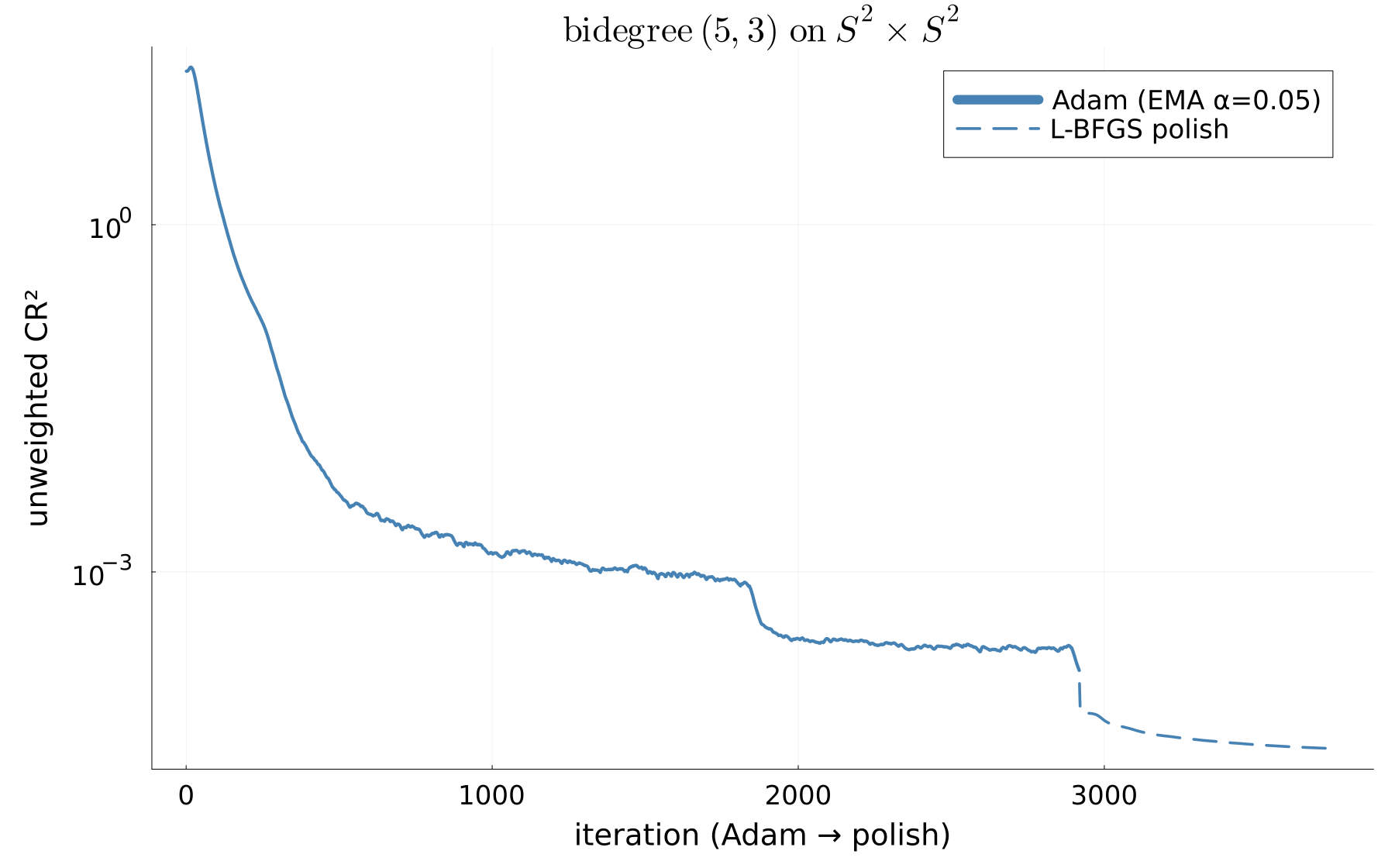}
    \caption{A graph of the training process.
    The first sharp cliff around iteration 1900 comes from to a halving of the learning rate after a 300-iteration plateau in $\mathcal L_{\mathrm{cr}}$.
    } 
    \label{fig:5_3_loss_plot}
\end{figure}

\section{Results for non-integrable $J$}

We next verify that our method can find $J$-holomorphic curves where $J$ is not integrable. Note here already we cannot compare to known solutions because there are no longer explicit solutions for nonintegrable $J$. Our only mechanism for checking whether we have found a solution to the $J$-holomorphic curve equation is checking that the Cauchy-Riemann loss is small.

In what follows, we will start with a known solution for an integrable $J_0$, then we consider what happens when we look at nearby nonintegrable $J$, and search for the corresponding $J$-holomorphic curves using our numerical schemes.

The symplectic manifold we will study is $S^2\times S^2$ and $J$-holomorphic curves in homology class $[S^2]\times[\mathrm{pt}]$ (bidegree $(0,1)$) and $[S^2\times \mathrm{pt}]+ [\mathrm{pt}\times S^2]$ (bidegree $(1,1)$). This symplectic manifold and the specific classes of $J$-holomorphic curves we study were featured prominently in symplectic geometry's history, for some examples see chapter 9 of \cite{McDuff_Salamon}.

We now specify the deformations of almost complex structure we will use. To leverage the fact $S^2\times S^2$ is embedded in $\mathbb{R}^6$, we will work mostly with vectors in $\mathbb{R}^6$ to simplify computations. For example, we take the standard almost complex structure $J_0$ to be actually a matrix acting on $\mathbb{R}^6$ that restricts to the standard almost complex structure in $S^2\times S^2$.

Let $Y_\epsilon$ be a 6 by 6 matrix, so that  $J_0Y_\epsilon=-J_0Y_\epsilon$. We choose $\epsilon$ small enough so that $1+1/2J_0Y_\epsilon$ is invertible.
We obtain an associated almost complex structure $J_\epsilon$ as
\begin{equation}
    J_{\epsilon} = (1+1/2J_0Y_\epsilon) J_0 (1+1/2J_0Y_\epsilon)^{-1}.\label{e:nonIntegrableJ}
\end{equation}
See the chapter 2 of the book \cite{wendl_lec} for this construction of deformation of almost complex structure.

In the case of the curves of bidegree $(1,1)$ we take 
\begin{equation}
    Y_\epsilon=\epsilon \left(\begin{matrix}
        A&B\\
        B^T&0
    \end{matrix}\right)
\end{equation}

The $A$ and $B$ are 3 by 3 matrices (depending smoothly on the domain $S^2\times S^2)$ that are built out of reflection in $\mathbb{R}^3$, parallel transport on $S^2$, and multiplication by smooth functions.

In the case curves with bidegree $(0,1)$ we take 
\begin{equation}
    Y_\epsilon=\epsilon \left(\begin{matrix}
        0&B\\
        B^T&C
    \end{matrix}\right)
\end{equation}
The $B$ and $C$ are also 3 by 3 matrices (depending smoothly on the domain $S^2\times S^2)$ that are built out of reflection in $\mathbb{R}^3$, parallel transport on $S^2$, and multiplication by smooth functions.

For detailed descriptions of the matrices $A$, $B$, $C$, refer to the functions \verb|_mixing_fibered_full| and \verb|_mixing_f2fibered_full| in \verb|src/targets/perturbed_product.jl| in the codebase; we note that because of the matrix inversion, some care had to be taken to ensure $J_\epsilon$ can be computed efficiently and differentiated through nicely.

We check numerically that $J_\epsilon$ preserves $T(S^2\times S^2)$, $J_\epsilon^2=-1$ on $T(S^2\times S^2)$, and that $J_\epsilon$ is tame. We verify $J_\epsilon$ is nonintegrable by computing the Nijenhuis tensor numerically.

To pin down specific $J$-holomorphic curves in an otherwise high dimensional moduli space we select marked points using the tetrahedron inscribed in the source $S^2$ with vertices
\begin{equation}
    v_1=\frac{1}{\sqrt 3}(1, 1, 1), \quad v_2=\frac{1}{\sqrt 3}(1, -1, -1), \quad v_3=\frac{1}{\sqrt 3}(-1, 1, -1), \quad v_4=\frac{1}{\sqrt 3}(-1, -1, 1).
\end{equation}

\subsection{The (0,1) foliation} \label{ss:gromovFoliation}

We fix an almost complex structure $J$ on $S^2\times S^2$ that corresponds to $\epsilon=0.9$ in the above description. We write the maps $u:S^2\rightarrow S^2\times S^2$ as $u=(u_1,u_2)$ where $u_i:S^2\rightarrow S^2$. Our marked points constraints are
\[
u_2(v_1)=v_1,\quad u_2(v_2)=v_2,\quad u_2(v_3)=v_3, \quad u_1(v_4)=c
\]
Then we consider a sequence of $J$-holomorphic curves corresponding to different values of $c \in\{c_0,c_1,..,c_6\}$. When $c=c_0=v_4$, by our choice of $J$, the $J$-holomorphic curve that satisfies the above marked point constraint is $u(z)=(v_4,z)$. This agrees with a holomorphic curve had we chosen the standard complex structure. Then we move the value of $c$ as a marked point constraint down the list along a meridian from $v_4$.
During this process we shall see the corresponding $J$-holomorphic curves pass through regions of $S^2\times S^2$ where the almost complex structure becomes more non-integrable. The different $J$-holomorphic curves corresponding to different values of $c$ are in fact leaves of a $J$-holomorphic foliation on $S^2\times S^2$.

The training and implementation of our algorithm proceeds as follows. We start with $c_0 = v_4$ and pretain our curve to the $J$-holomorphic map $u(z) = (c_0,z)$. Then we take then sequence of $c\in \{c_0,c_1,..,c_6\}$ where the points move along the meridian. We look for a sequence $u^i$ of $J$-holomorphic maps satisfying the same constraints at $v_1,v_2,v_3$ as before and satisfies $u^i_1(v_4)=c_i$ . At the $i$-th step, we assume we have pretrained to the $u^i$ map, and run our numerical algorithm to find the $i+1$-th map (i.e. we move the marked point from $c_i$ to $c_{i+1}$ in the loss function, and search for a new map that is $J$-holomorphic and satisfies this new marked point constraint). Because successive curves are close to each other and have nearby marked points, for this problem we use paradigm 1 with a depth 2 neural network.

We plot the results for each of $c_0$, $c_2$, and $c_5$ as follows.
We put spherical coordinates on both the domain $S^2$ ($\phi,\theta$) and target $S^2\times S^2$ $(\phi_1,\theta_1,\phi_2,\theta_2)$. This means to specify a given map we need 4 images. We then plot the Cauchy Riemann loss that measures how $J$-holomorphic the map is. 
Finally we measure the distance (using the distance in $\mathbb{R}^6$) of the map we found against $v^i(z) = (c_i,z)$  as a measure of how much changing the almost complex structure to be non-integrable has deformed the holomorphic maps. 

We also include a table for the results for $c_1,c_2,...,c_6$ indicating Cauchy-Riemann loss and marked point loss.

\begin{table}[H]
    \centering
    \begin{tabular}{c|c|c|c|c|c|c}
    leaf & $\theta$ & mean $\mathcal{L}_{\mathrm{cr}}$ & $\mathcal{L}_{\mathrm{mp}}$ & mean $\operatorname{dist}(u,(c_i,z))$
         & factor 1 deg & factor 2 deg \\
         \hline
$1$ & $0$ & $9.1877 \times 10^{-7}$ & $4.4974 \times 10^{-11}$ & $4.4361 \times 10^{-5}$ & $0.0$ & $1.0$\\
$2$ & $\pi/10$ & $3.16 \times 10^{-7}$ & $3.999 \times 10^{-12}$ & $0.050174$ & $0.0$ & $0.99995$\\
$3$ & $\pi/5$ & $3.1142 \times 10^{-7}$ & $2.2051 \times 10^{-11}$ & $0.20115$ & $1.0 \times 10^{-5}$ & $0.99981$\\
$4$ & $3\pi/10$ & $4.5891 \times 10^{-7}$ & $1.6001 \times 10^{-11}$ & $0.38501$ & $2.0 \times 10^{-5}$ & $0.99969$\\
$5$ & $2\pi/5$ & $6.284 \times 10^{-7}$ & $2.0715 \times 10^{-11}$ & $0.49994$ & $3.0 \times 10^{-5}$ & $0.99967$\\
$6$ & $\pi/2$ & $5.7931 \times 10^{-7}$ & $1.2116 \times 10^{-11}$ & $0.52685$ & $3.0 \times 10^{-5}$ & $0.99969$
    \end{tabular}
    \caption{Results for the six leaves of the foliation.}
    \label{tab:01foliation}
\end{table}

\begin{figure}[H]
    \centering
    \includegraphics[width=0.8\linewidth]{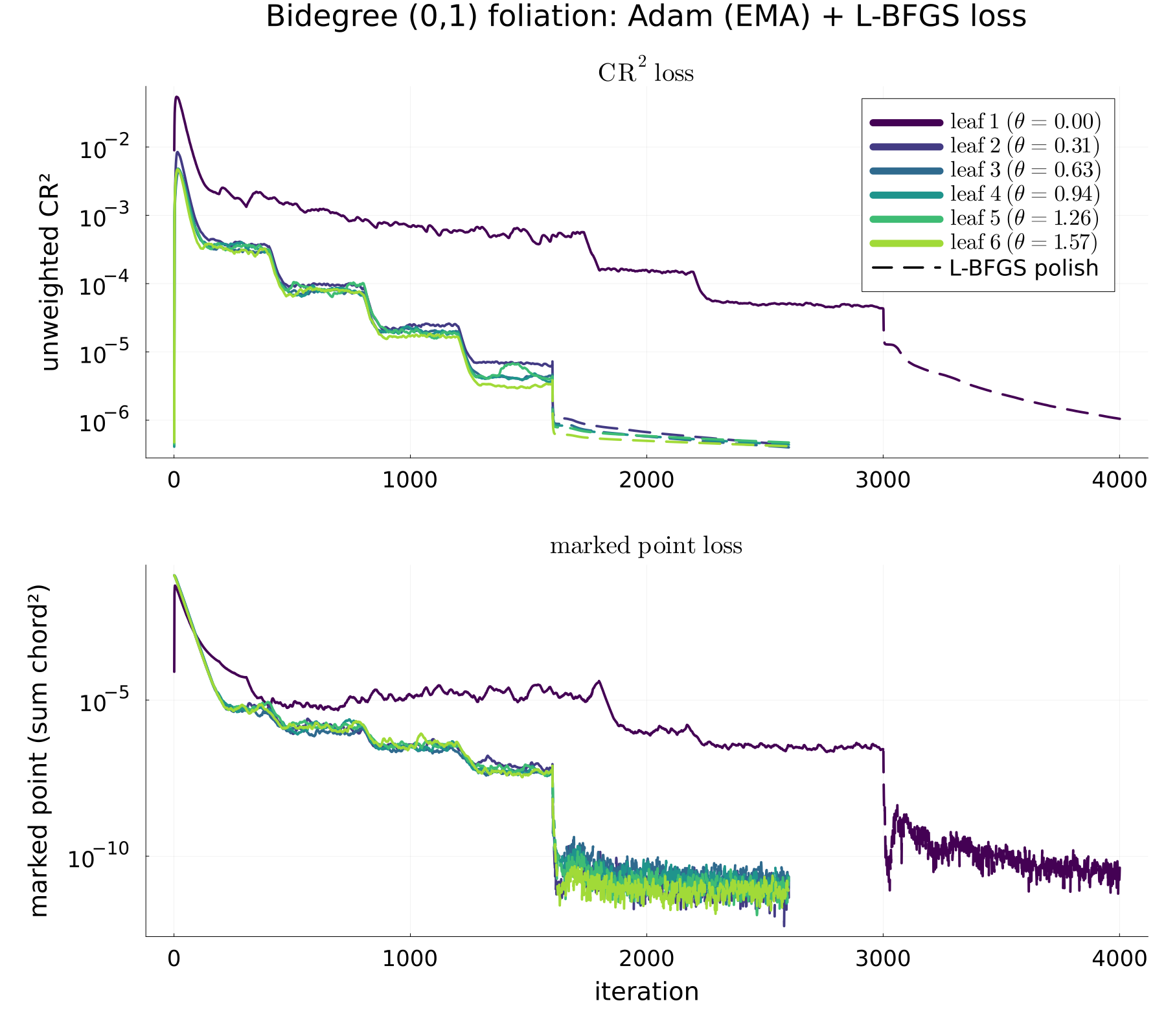}
    \caption{Loss curves for $\mathcal{L}_{\mathrm{cr}}$ and $\mathcal{L}_{\mathrm{mp}}$ for the training of all six leafs.
    Note that the initial leaf trains for the full $3000$ allotted Adam iterations;  the initial leaf is pretrained to $(c_0,z)$ only in $L^2$ and needs the full time to become holomorphic.
    All leaves after the first benefit from starting from the previous leaf (and enter the L-BFGS polish early after Adam plateaus).
    }
    \label{fig:01foliation_loss}
\end{figure}

\subsection{The $(1,1)$ $J_\epsilon$-holomorphic curves}\label{ss:bifurcation}

We consider $J$-holomorphic curves $u:S^2\rightarrow S^2\times S^2$ of bidegree $(1,1)$. To cut down the dimension of the moduli space we impose marked point constraints
\[
u(v_1) = (v_1,v_1),\quad u(v_2)=(v_2,v_2), \quad u(v_3)= (v_3,v_3)
\]

When $\epsilon=0$, the almost complex structure $J_\epsilon$ given by~\eqref{e:nonIntegrableJ} reduces to the standard product almost complex structure $J_0$ on $S^2\times S^2$.
The diagonal map is the unique $J_0$-holomorphic curve with bidegree $(1,1)$ satisfying these marked point. We start with an $J_0$-holomorphic curve satisfying the marked point constraints and watch it deform as we turn up $\epsilon$ and enforce the same marked point constraints. We consider  $\epsilon_i\in [0,0.3, 0.6, 0.9, 1.2]$.

Our algorithm is implemented as follows. We first learn the holomorphic ground truth $u^0=(z,z)$ when $\epsilon=0$; this is done by following paradigm 1 above, first pretraining to $(z,z)$ in the $L^2$ norm and then requiring that it be holomorphic.

For the $i$th stage, we start with the curve $u^i$ which is $J_{\epsilon_i}$-holomorphic satisfying the marked point constraints, then perturb the almost complex structure to $J_{\epsilon_{i+1}}$. 
We then use paradigm 1 with  a depth 2 network to find the $J_{\epsilon_{i+1}}$-holomorphic curve $u^{i+1}$ that satisfies the same marked point constraints.
Because the marked point constraints are satisfied exactly for each $\epsilon$ and the initial pretrain target $(z,z)$ also satisfies them, we do not need to use the broader exploratory power of paradigm 2.

We plot the results below. We use spherical coordinates on the domain $S^2$, where each point is specified by two angles $(\phi,\theta)$. We likewise use spherical coordinates on taget $S^2\times S^2$ with angles $(\phi_1,\theta_1,\phi_2,\theta_2)$. Hence to specify each map we need 4 panels. We also show the Cauchy-Riemann loss, and the distance of our $J_\epsilon$-holomorphic curve $u$ from the diagonal map $z\rightarrow (z,z)$ in $\mathbb{R}^6$ as a way to show how much the curve has deformed as a consequence of changing $J_0$.

\begin{table}[H]
    \centering
    \begin{tabular}{c|c|c|c|c|c}
    $\epsilon$ &   $\mathcal{L}_{\mathrm{cr}}$ & $\mathcal{L}_{\mathrm{mp}}$ & mean $\operatorname{dist}(u,(z,z))$
         & factor 1 deg & factor 2 deg \\
         \hline
$0.0$ & $1.04 \times 10^{-6}$ & $5.5731 \times 10^{-11}$ & $3.74 \times 10^{-5}$ & $1.0$ & $1.0$\\
$0.3$ & $4.13 \times 10^{-7}$ & $9.8751 \times 10^{-12}$ & $0.191$ & $0.99985$ & $1.0$\\
$0.6$ & $7.99 \times 10^{-7}$ & $1.7803 \times 10^{-12}$ & $0.377$ & $0.99977$ & $1.0$\\
$0.9$ & $4.82 \times 10^{-6}$ & $1.2207 \times 10^{-10}$ & $0.561$ & $0.99976$ & $0.99999$\\
$1.2$ & $\approx 0.5$ & $6.7679 \times 10^{-8}$ & $0.746$ & $0.9999$ & $1.0007$
    \end{tabular}
    \caption{Statistics for the trained network for each of the five values of $\epsilon$ used in the continuation; the Cauchy-Riemann loss is computed on a more refined grid that is held out during training.
    For $\epsilon=1.2$, the dependence of $\mathcal{L}_{\mathrm{cr}}$ on the grid is stronger and we report fewer digits as a result.
    }
    \label{tab:11continuation}
\end{table}

We note $\epsilon=1.2$ yielded substantially worse results. 
Its Cauchy-Riemann loss is at a value where we usually declare training has failed. 
We expect this is because $\epsilon=1.2$ is near the value where $1+1/2 J_0Y_\epsilon$ is no longer invertible, as measured by the presence of very small eigenvalues/singular values. 
This means we expect very large numerical entries in the $J_\epsilon$ matrix, which would require substantially smaller learning rates for the numerical algorithm to yield good results.

\begin{figure}[H]
    \centering
    \includegraphics[width=0.8\textwidth]{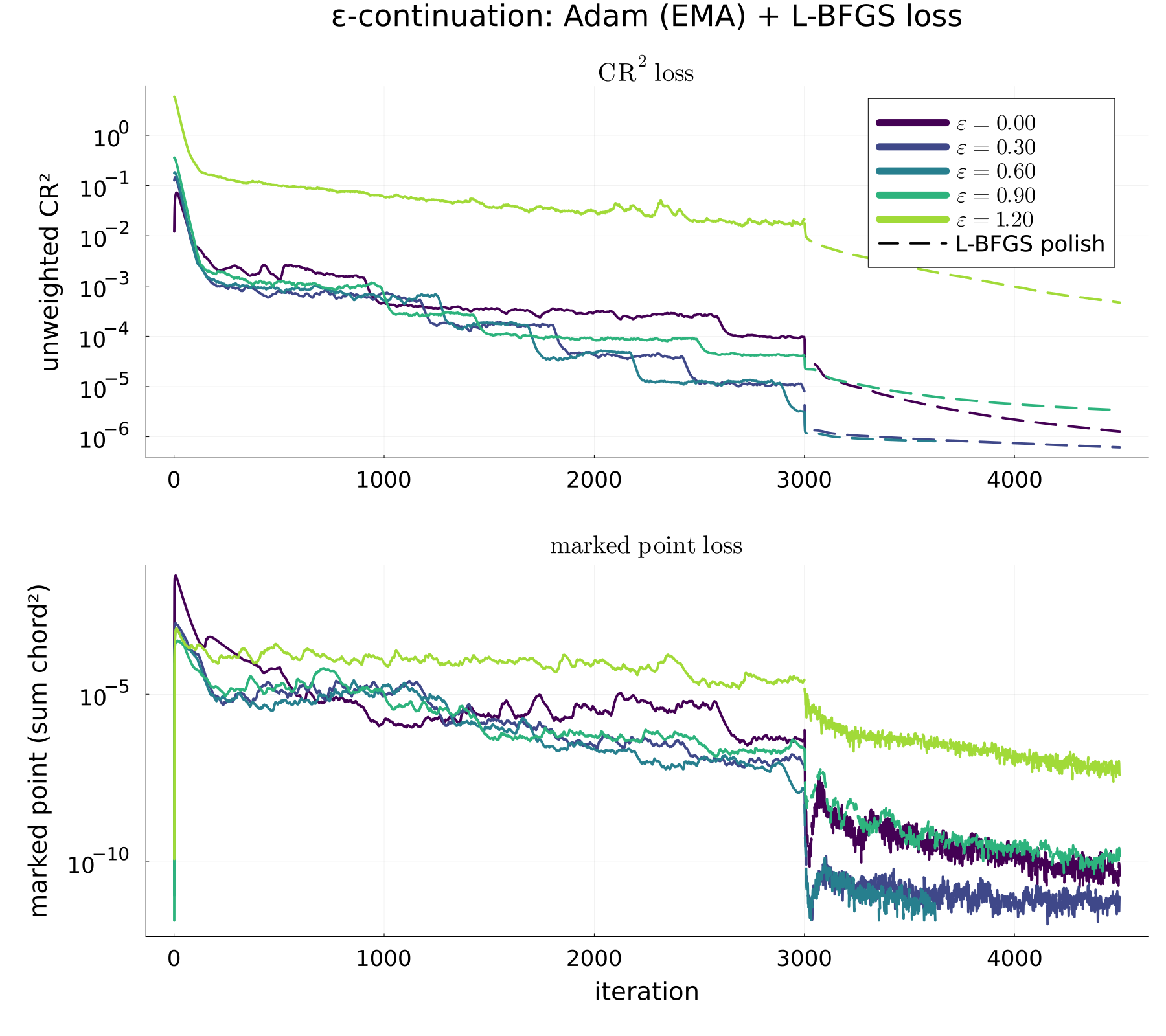}
    \caption{Loss curves for $\mathcal{L}_{\mathrm{cr}}$ and $\mathcal{L}_{\mathrm{mp}}$ for the training of all five values of $\epsilon$.
    The comparatively high losses for $\epsilon=0$ are due to the fact that later iterations benefit from starting from the previous curve which is nearly $J_{\epsilon_{i-1}}$-holomorphic, while pretraining to the diagonal identity is only done in $L^2$.
    }
    \label{fig:continuation_loss}
\end{figure}

\section{Future directions}
\label{sec:future}
The subject of symplectic geometry has developed greatly since the initial introduction of $J$-holomorphic curves by Gromov. In particular people have introduced methods to study $J$-holomorphic curves with boundary and punctures, and counts of these curves assemble into algebraic gadgets like homology theories and categories. Developments in the community led to rather sophisticated tools to extract geometric information from these more involved invariants. Consequently there are several natural extensions of this paper:
\begin{enumerate}
\item Build a $J$-holomorphic curve finder for curves with interior punctures, for example those arising in Hamiltonian Floer homology \cite{audin_damian} or symplectic field theory (we refer to the textbook \cite{sft}), with boundary conditions on Lagrangians (open Gromov-Witten invariants \cite{open_gw}) and boundary punctures (Lagrangian Floer homology See \cite{auroux2013}). 

\item To extract the most interesting information from $J$-holomorphic curves, it is usually not enough to simply show these curves exist, it is often important to explicitly count how many $J$-holomorphic curves there are.  Our results in Section \ref{sec:wp} suggest that this is very much achievable with the tools we have built.

\item The possible degeneration of $J$-holomorphic curves is understood via Gromov's compactification (see \cite{McDuff_Salamon}), it would be interesting to see degeneration of $J$-holomorphic curves explicitly in numerical simulations.

\item In the earlier days of the development of symplectic geometry, variational techniques saw widespread usage in order to find periodic orbits of (Hamiltonian) vector fields on symplectic manifolds \cite{hofer_zehnder}. Today variational techniques see less usage compared to the Floer theory package. However, the results of this paper are a demonstration that variational techniques with enough computing power do work when tasked with finding interesting objects in symplectic geometry, in this case $J$-holomorphic curves. It would be interesting to revisit on the pure math side the usage of variational techniques to find $J$-holomorphic curves.

\end{enumerate}

\section*{Acknowledgments}
The code was written with the help of Claude (primarily the Opus 4.8 and Opus 5 models). 
We read the code line-by-line and are responsible for its correctness. 
The observation that $\wp^{(1)}$ and $\wp^{(2)}$ are the unique two degree 2 holomorphic maps $T^2\rightarrow S^2$ satisfying the marked point constraints was originally suggested by Claude/Gemini. 
We verified this statement by hand. Everything else in this project was done by the two authors.

We would like to thank Kyler Siegel for helpful conversations, Elliot Kienzle for comments on an earlier draft of the paper, and Javier Gomez-Serrano for suggesting a more detailed comparison to non-PINN methods. We are inspired by the creators of the symplectic embedding machine, Alex Gajewski, Eli Goldin, Semon Rezchikov, Jakwanul Safin, Navtej Singh, Kyler Siegel and Junhui Zhang.

\bibliographystyle{plain}
\bibliography{jholomorphic_pinns}

@article {raissi2019PINNs,
    AUTHOR = {Raissi, M. and Perdikaris, P. and Karniadakis, G. E.},
     TITLE = {Physics-informed neural networks: a deep learning framework
              for solving forward and inverse problems involving nonlinear
              partial differential equations},
   JOURNAL = {J. Comput. Phys.},
  FJOURNAL = {Journal of Computational Physics},
    VOLUME = {378},
      YEAR = {2019},
     PAGES = {686--707},
      ISSN = {0021-9991,1090-2716},
   MRCLASS = {65M70 (68T05)},
  MRNUMBER = {3881695},
       DOI = {10.1016/j.jcp.2018.10.045},
       URL = {https://doi.org/10.1016/j.jcp.2018.10.045},
}

@inproceedings{tancik2020fourierFeatures,
  author       = {Matthew Tancik and
                  Pratul P. Srinivasan and
                  Ben Mildenhall and
                  Sara Fridovich{-}Keil and
                  Nithin Raghavan and
                  Utkarsh Singhal and
                  Ravi Ramamoorthi and
                  Jonathan T. Barron and
                  Ren Ng},
  editor       = {Hugo Larochelle and
                  Marc'Aurelio Ranzato and
                  Raia Hadsell and
                  Maria{-}Florina Balcan and
                  Hsuan{-}Tien Lin},
  title        = {Fourier Features Let Networks Learn High Frequency Functions in Low
                  Dimensional Domains},
  booktitle    = {Advances in Neural Information Processing Systems 33: Annual Conference
                  on Neural Information Processing Systems 2020, NeurIPS 2020, December
                  6-12, 2020, virtual},
  year         = {2020},
  url          = {https://proceedings.neurips.cc/paper/2020/hash/55053683268957697aa39fba6f231c68-Abstract.html},
  bibsource    = {dblp computer science bibliography, https://dblp.org}
}

@inproceedings{rahaman2019spectralBias,
  author       = {Nasim Rahaman and
                  Aristide Baratin and
                  Devansh Arpit and
                  Felix Draxler and
                  Min Lin and
                  Fred A. Hamprecht and
                  Yoshua Bengio and
                  Aaron C. Courville},
  editor       = {Kamalika Chaudhuri and
                  Ruslan Salakhutdinov},
  title        = {On the Spectral Bias of Neural Networks},
  booktitle    = {Proceedings of the 36th International Conference on Machine Learning,
                  {ICML} 2019, 9-15 June 2019, Long Beach, California, {USA}},
  series       = {Proceedings of Machine Learning Research},
  volume       = {97},
  pages        = {5301--5310},
  publisher    = {{PMLR}},
  year         = {2019},
  url          = {http://proceedings.mlr.press/v97/rahaman19a.html},
  bibsource    = {dblp computer science bibliography, https://dblp.org}
}

@article{hendrycks2016gelu,
  title={Gaussian error linear units (GELUs)},
  author={Hendrycks, Dan and Gimpel, Kevin},
  journal={arXiv preprint arXiv:1606.08415},
  year={2016}
}

@inproceedings{diederik2015adam,
  author       = {Diederik P. Kingma and
                  Jimmy Ba},
  editor       = {Yoshua Bengio and
                  Yann LeCun},
  title        = {Adam: {A} Method for Stochastic Optimization},
  booktitle    = {3rd International Conference on Learning Representations, {ICLR} 2015,
                  San Diego, CA, USA, May 7-9, 2015, Conference Track Proceedings},
  year         = {2015},
  url          = {http://arxiv.org/abs/1412.6980},
  bibsource    = {dblp computer science bibliography, https://dblp.org}
}

@article {liu1989lbfgs,
    AUTHOR = {Liu, Dong C. and Nocedal, Jorge},
     TITLE = {On the limited memory {BFGS} method for large scale
              optimization},
   JOURNAL = {Math. Programming},
  FJOURNAL = {Mathematical Programming},
    VOLUME = {45},
      YEAR = {1989},
    NUMBER = {3},
     PAGES = {503--528},
      ISSN = {0025-5610,1436-4646},
   MRCLASS = {90C30 (65K05)},
  MRNUMBER = {1038245},
       DOI = {10.1007/BF01589116},
       URL = {https://doi.org/10.1007/BF01589116},
}

@book {McDuff_Salamon,
    AUTHOR = {McDuff, Dusa and Salamon, Dietmar},
     TITLE = {{$J$}-holomorphic curves and symplectic topology},
    SERIES = {American Mathematical Society Colloquium Publications},
    VOLUME = {52},
   EDITION = {Second},
 PUBLISHER = {American Mathematical Society, Providence, RI},
      YEAR = {2012},
     PAGES = {xiv+726},
      ISBN = {978-0-8218-8746-2},
   MRCLASS = {53D45 (32Q65 53D35)},
  MRNUMBER = {2954391},
MRREVIEWER = {Mark\ Alan\ Branson},
}

@article {Gromov,
    AUTHOR = {Gromov, M.},
     TITLE = {Pseudo holomorphic curves in symplectic manifolds},
   JOURNAL = {Invent. Math.},
  FJOURNAL = {Inventiones Mathematicae},
    VOLUME = {82},
      YEAR = {1985},
    NUMBER = {2},
     PAGES = {307--347},
    
}

@book {intro_sympl,
    AUTHOR = {McDuff, Dusa and Salamon, Dietmar},
     TITLE = {Introduction to symplectic topology},
    SERIES = {Oxford Graduate Texts in Mathematics},
   EDITION = {Third},
 PUBLISHER = {Oxford University Press, Oxford},
      YEAR = {2017},
     PAGES = {xi+623},
      ISBN = {978-0-19-879490-5; 978-0-19-879489-9},
   MRCLASS = {53D35 (53D40 57R17 57R57 57R58)},
  MRNUMBER = {3674984},
MRREVIEWER = {Hansj\"org\ Geiges},
       DOI = {10.1093/oso/9780198794899.001.0001},
       URL = {https://doi.org/10.1093/oso/9780198794899.001.0001},
}

@article {wendl,
    AUTHOR = {Wendl, Chris},
     TITLE = {Automatic transversality and orbifolds of punctured
              holomorphic curves in dimension four},
   JOURNAL = {Comment. Math. Helv.},
  FJOURNAL = {Commentarii Mathematici Helvetici. A Journal of the Swiss
              Mathematical Society},
    VOLUME = {85},
      YEAR = {2010},
    NUMBER = {2},
     PAGES = {347--407},
      ISSN = {0010-2571,1420-8946},
   MRCLASS = {32Q65 (53D45 57R17)},
  MRNUMBER = {2595183},
MRREVIEWER = {Umberto\ Leone\ Hryniewicz},
       DOI = {10.4171/CMH/199},
       URL = {https://doi.org/10.4171/CMH/199},
}

@book {audin_damian,
    AUTHOR = {Audin, Mich\`ele and Damian, Mihai},
     TITLE = {Morse theory and {F}loer homology},
    SERIES = {Universitext},
      NOTE = {Translated from the 2010 French original by Reinie Ern\'e},
 PUBLISHER = {Springer, London; EDP Sciences, Les Ulis},
      YEAR = {2014},
     PAGES = {xiv+596},
      ISBN = {978-1-4471-5495-2; 978-1-4471-5496-9; 978-2-7598-0704-8},
   MRCLASS = {53-02 (53D40 58E05)},
  MRNUMBER = {3155456},
MRREVIEWER = {Sonja\ Hohloch},
       DOI = {10.1007/978-1-4471-5496-9},
       URL = {https://doi.org/10.1007/978-1-4471-5496-9},
}

@misc{sft,
      title={Lectures on Symplectic Field Theory}, 
      author={Chris Wendl},
      year={2016},
      eprint={1612.01009},
      archivePrefix={arXiv},
      primaryClass={math.SG},
      url={https://arxiv.org/abs/1612.01009}, 
}

@misc{wendl_lec, author={Chirs Wendl}, title={Lectures on Holomorphic Curves in
Symplectic and Contact Geometry},year={2015}, url={https://www.mathematik.hu-berlin.de/~wendl/pub/jhol_bookv33.pdf}}

@misc{auroux2013,
      title={A beginner's introduction to Fukaya categories}, 
      author={Denis Auroux},
      year={2013},
      eprint={1301.7056},
      archivePrefix={arXiv},
      primaryClass={math.SG},
      url={https://arxiv.org/abs/1301.7056}, 
}

@book {hofer_zehnder,
    AUTHOR = {Hofer, Helmut and Zehnder, Eduard},
     TITLE = {Symplectic invariants and {H}amiltonian dynamics},
    SERIES = {Birkh\"auser Advanced Texts: Basler Lehrb\"ucher.
              [Birkh\"auser Advanced Texts: Basel Textbooks]},
 PUBLISHER = {Birkh\"auser Verlag, Basel},
      YEAR = {1994},
     PAGES = {xiv+341},
      ISBN = {3-7643-5066-0},
   MRCLASS = {58-02 (34C25 57R15 58E05 58F05 70H05)},
  MRNUMBER = {1306732},
MRREVIEWER = {Daniel\ M.\ Burns, Jr.},
       DOI = {10.1007/978-3-0348-8540-9},
       URL = {https://doi.org/10.1007/978-3-0348-8540-9},
}

@book {arnold,
    AUTHOR = {Arnold, V. I.},
     TITLE = {Mathematical methods of classical mechanics},
    SERIES = {Graduate Texts in Mathematics},
    VOLUME = {60},
      NOTE = {Translated from the Russian by K. Vogtmann and A. Weinstein},
 PUBLISHER = {Springer-Verlag, New York-Heidelberg},
      YEAR = {1978},
     PAGES = {x+462},
      ISBN = {0-387-90314-3},
   MRCLASS = {58F05 (70.58)},
  MRNUMBER = {690288},
MRREVIEWER = {J.\ S.\ Joel},
}

@book {mirror_symmetry,
    AUTHOR = {Hori, Kentaro and Katz, Sheldon and Klemm, Albrecht and
              Pandharipande, Rahul and Thomas, Richard and Vafa, Cumrun and
              Vakil, Ravi and Zaslow, Eric},
     TITLE = {Mirror symmetry},
    SERIES = {Clay Mathematics Monographs},
    VOLUME = {1},
      NOTE = {With a preface by Vafa},
 PUBLISHER = {American Mathematical Society, Providence, RI; Clay
              Mathematics Institute, Cambridge, MA},
      YEAR = {2003},
     PAGES = {xx+929},
      ISBN = {0-8218-2955-6},
   MRCLASS = {14J32 (14N35 32Q25 81T30 81T60)},
  MRNUMBER = {2003030},
MRREVIEWER = {Marcos\ Mari\~no},
}

@article {open_gw,
    AUTHOR = {Solomon, Jake P. and Tukachinsky, Sara B.},
     TITLE = {Point-like bounding chains in open {G}romov-{W}itten theory},
   JOURNAL = {Geom. Funct. Anal.},
  FJOURNAL = {Geometric and Functional Analysis},
    VOLUME = {31},
      YEAR = {2021},
    NUMBER = {5},
     PAGES = {1245--1320},
      ISSN = {1016-443X,1420-8970},
   MRCLASS = {53D45 (14N10 14N35 53D12 53D37)},
  MRNUMBER = {4356703},
MRREVIEWER = {Cheng-Yong\ Du},
       DOI = {10.1007/s00039-021-00583-3},
       URL = {https://doi.org/10.1007/s00039-021-00583-3},
}

@misc{elliot,
  author = {Elliot Kienzle},
 
  howpublished = {\url{https://chessapig.github.io/floer}},

}

@article{min_surface_survey,
    author = {Mifodijus Sapagovas and Vytautas Būda and Saulius Maskeliūnas and Olga Štikonienė and Artūras Štikonas},
    title = {Minimal Surfaces and the Plateau Problem: Numerical Methods and Applications},
    journal = {Informatica},
    volume = {35},
    number = {2},
    year = {2024},
    pages = {401--420},
    doi = {10.15388/24-INFOR552},
    issn = {0868-4952},
    publisher = {Vilnius University Institute of Data Science and Digital Technologies}
}

@article{Pinn_ricci,
  author       = {Aarjav Jain and
                  Challenger Mishra and
                  Pietro Li{\`{o}}},
  title        = {A physics-informed search for metric solutions to Ricci flow, their
                  embeddings, and visualisation},
  journal      = {CoRR},
  volume       = {abs/2212.05892},
  year         = {2022},
  url          = {https://doi.org/10.48550/arXiv.2212.05892},
  doi          = {10.48550/ARXIV.2212.05892},
  eprinttype   = {arXiv},
  eprint       = {2212.05892},
  bibsource    = {dblp computer science bibliography, https://dblp.org}
}

@misc{2026usersguidepinnsgeometric,
      title={A user's guide to PINNs in geometric analysis: lessons from the asymptotic Plateau problem}, 
      author={Tancredi Schettini Gherardini},
      year={2026},
      eprint={2607.28733},
      archivePrefix={arXiv},
      primaryClass={math.DG},
      url={https://arxiv.org/abs/2607.28733}, 
}

@misc{2026minimalsurfacesknotsneural,
      title={Minimal surfaces, Knots, and Neural Networks}, 
      author={Tancredi Schettini Gherardini and Marco Usula},
      year={2026},
      eprint={2605.26234},
      archivePrefix={arXiv},
      primaryClass={math.DG},
      url={https://arxiv.org/abs/2605.26234}, 
}

@misc{corts2026machinelearningapproachnirenberg,
      title={A Machine Learning Approach to the Nirenberg Problem}, 
      author={Gianfranco Cortés and Maria Esteban-Casadevall and Yueqing Feng and Jonas Henkel and Edward Hirst and Tancredi Schettini Gherardini and Alexander G. Stapleton},
      year={2026},
      eprint={2602.12368},
      archivePrefix={arXiv},
      primaryClass={cs.LG},
      url={https://arxiv.org/abs/2602.12368}, 
}

@misc{fang2021,
      title={A Physics-Informed Neural Network Framework For Partial Differential Equations on 3D Surfaces: Time-Dependent Problems}, 
      author={Zhiwei Fang and Justin Zhang and Xiu Yang},
      year={2021},
      eprint={2103.13878},
      archivePrefix={arXiv},
      primaryClass={cs.LG},
      url={https://arxiv.org/abs/2103.13878}, 
}

@InProceedings{douglas22a,
  title = 	 {Numerical Calabi-Yau metrics from holomorphic networks},
  author =       {Douglas, Michael and Lakshminarasimhan, Subramanian and Qi, Yidi},
  booktitle = 	 {Proceedings of the 2nd Mathematical and Scientific Machine Learning Conference},
  pages = 	 {223--252},
  year = 	 {2022},
  editor = 	 {Bruna, Joan and Hesthaven, Jan and Zdeborova, Lenka},
  volume = 	 {145},
  series = 	 {Proceedings of Machine Learning Research},
  month = 	 {16--19 Aug},
  publisher =    {PMLR},
  url = 	 {https://proceedings.mlr.press/v145/douglas22a.html}
}

@article {MR4608987,
    AUTHOR = {Wang, Y. and Lai, C.-Y. and G\'omez-Serrano, J. and
              Buckmaster, T.},
     TITLE = {Asymptotic self-similar blow-up profile for three-dimensional
              axisymmetric {E}uler equations using neural networks},
   JOURNAL = {Phys. Rev. Lett.},
  FJOURNAL = {Physical Review Letters},
    VOLUME = {130},
      YEAR = {2023},
    NUMBER = {24},
     PAGES = {Paper No. 244002, 6},
      ISSN = {0031-9007,1079-7114},
   MRCLASS = {76B03},
  MRNUMBER = {4608987},
       DOI = {10.1103/physrevlett.130.244002},
       URL = {https://doi.org/10.1103/physrevlett.130.244002},
}

@article {MR4761208,
    AUTHOR = {Berglund, Per and Butbaia, Giorgi and H\"ubsch, Tristan and
              Jejjala, Vishnu and Mayorga Pe\~na, Dami\'an and Mishra,
              Challenger and Tan, Justin},
     TITLE = {Machine learned {C}alabi-{Y}au metrics and curvature},
   JOURNAL = {Adv. Theor. Math. Phys.},
  FJOURNAL = {Advances in Theoretical and Mathematical Physics},
    VOLUME = {27},
      YEAR = {2023},
    NUMBER = {4},
     PAGES = {1107--1158},
      ISSN = {1095-0761,1095-0753},
   MRCLASS = {32Q25 (14J32 32G05 53C55 68T07)},
  MRNUMBER = {4761208},
MRREVIEWER = {Nadaniela\ Egidi},
       DOI = {10.4310/atmp.2023.v27.n4.a3},
       URL = {https://doi.org/10.4310/atmp.2023.v27.n4.a3},
}

@article{Gerdes_2023,
doi = {10.1088/2632-2153/acdc84},
url = {https://doi.org/10.1088/2632-2153/acdc84},
year = {2023},
month = {jun},
publisher = {IOP Publishing},
volume = {4},
number = {2},
pages = {025031},
author = {Gerdes, Mathis and Krippendorf, Sven},
title = {CYJAX: A package for Calabi-Yau metrics with JAX},
journal = {Machine Learning: Science and Technology}
}

@article {deepritz,
    AUTHOR = {Berglund, Per and Butbaia, Giorgi and H\"ubsch, Tristan and
              Jejjala, Vishnu and Mayorga Pe\~na, Dami\'an and Mishra,
              Challenger and Tan, Justin},
     TITLE = {The Deep Ritz Method: A Deep Learning-Based Numerical Algorithm for Solving Variational Problems},
   JOURNAL = {Communications in Mathematics and Statistics},
 
    VOLUME = {6},
      YEAR = {2018},
    NUMBER = {1},
     PAGES = {1--12}
      
}

@article {Levenberg,
    AUTHOR = {Levenberg, Kenneth},
     TITLE = {A method for the solution of certain non-linear problems in
              least squares},
   JOURNAL = {Quart. Appl. Math.},
  FJOURNAL = {Quarterly of Applied Mathematics},
    VOLUME = {2},
      YEAR = {1944},
     PAGES = {164--168},
      ISSN = {0033-569X,1552-4485},
   MRCLASS = {65.0X},
  MRNUMBER = {10666},
MRREVIEWER = {T.\ E.\ Sterne},
       DOI = {10.1090/qam/10666},
       URL = {https://doi.org/10.1090/qam/10666},
}

@article {Marquardt,
    AUTHOR = {Marquardt, Donald W.},
     TITLE = {An algorithm for least-squares estimation of nonlinear
              parameters},
   JOURNAL = {J. Soc. Indust. Appl. Math.},
  FJOURNAL = {Journal of the Society for Industrial and Applied Mathematics},
    VOLUME = {11},
      YEAR = {1963},
     PAGES = {431--441},
      ISSN = {0368-4245},
   MRCLASS = {62.20},
  MRNUMBER = {153071},
MRREVIEWER = {M.\ Atiqullah},
}
    
\end{document}